\documentclass[ascmac,11pt,twoside]{interact}
\usepackage{epstopdf}
\usepackage[caption=false]{subfig}
\usepackage[natbibapa,nodoi]{apacite}
\theoremstyle{plain}
\newtheorem{theorem}{Theorem}[section]
\newtheorem{lemma}[theorem]{Lemma}
\newtheorem{corollary}[theorem]{Corollary}

\theoremstyle{definition}
\newtheorem{definition}[theorem]{Definition}

\theoremstyle{remark}

\usepackage{graphicx}
\usepackage{amsmath,amssymb,amsfonts}
\usepackage{amsthm}
\usepackage{color}
\usepackage{dcolumn}
\newcolumntype{d}[1]{D{.}{\cdot}{#1}}
\newcolumntype{.}{D{.}{.}{-1}}
\newcolumntype{,}{D{,}{,}{-1}}

\usepackage{hhline}
\usepackage{multirow}

\theoremstyle{remark}
\theoremstyle{definition}

\newcommand{\bb}[1]{\mathbb{#1}}  

\newfont{\bg}{cmr10 scaled\magstep4} 

\newcommand{\bigzerou}{\smash{\lower1.7ex\hbox{\bg 0}}}
\usepackage{bm}
\usepackage{upgreek}
\usepackage{longtable}
\usepackage{here}
\usepackage{color}
\usepackage{lineno}
\usepackage{blkarray, bigstrut}
\usepackage{mathrsfs}
\newcommand{\bbR}{\mathbb{R}}

\newcommand{\Pexfa}{ P_{\scalebox{0.5}{\it EXFA}} }
\newcommand{\Pcon}{Q_{\scalebox{0.5}{\it CON}}}
\usepackage[%
 setpagesize=false,%
 bookmarks=true,%
 bookmarksdepth=tocdepth,%
 bookmarksnumbered=true,%
 colorlinks=true,%
 pdftitle={},%
 pdfsubject={},%
 pdfauthor={},%
 pdfkeywords={}%
]{hyperref}
\hypersetup{
    citecolor=blue,
    linkcolor=blue,
}

\begin{document}


\title{Closest targets in Russell graph measure of strongly monotonic efficiency for an extended facet production possibility set}

\author{
\name{Kazuyuki Sekitani\textsuperscript{a}\thanks{K. Sekitani (kazuyuki-sekitani@st.seikei.ac.jp)} and Yu Zhao\textsuperscript{b}\thanks{Y. Zhao (yu.zhao@rs.tus.ac.jp)}}
\affil{\textsuperscript{a}Seikei University, 3-3-1 Kichijoji-Kitamachi, Musashino-shi, Tokyo 180-8633, Japan;\\ \textsuperscript{b}Tokyo University of Science, 1-3 Kagurazaka, Shinjuku-ku, Tokyo 162-8601, Japan}
}

\maketitle

\begin{abstract}
This study modifies the Russell graph measure of efficiency to find the closest target from the assessed decision-making unit to the efficient frontier of an empirical production possibility set using an extended facet approach. 
The closest target requires a single improvement in either an output term or an input term, and needs some  positive  levels of all inputs to achieve the output target.
The practical advantage of the proposed efficiency measure is demonstrated numerically in comparison to other existing non-radial efficiency measures.
\end{abstract}



\section{Introduction}
Data envelopment analysis (DEA) is a mathematical programming method to measure the relative efficiency of the production activities called decision-making units (DMUs). 
Some objective functions of DEA models are nonlinear with respect to the levels of input and output improvement.
For example,  the DEA model of the Russell graph measure~\citep[RM;][]{fare1985measurement}{}{}  has a sum of linear ratios as the objective function and 
the slacks-based measure~\citep[SBM;][]{tone2001slacks}{}{} is formulated as a DEA model with  a linear fractional function. 
SBM  is also known as the enhanced Russell graph measure~\citep[ERGM;][]{pastor1999enhanced}{}{}.
Recently,~\citet{salahi2019robust} developed  robust RM and ERGM DEA models under certain types of  uncertainties, while~\citet{TONE2020560} proposed a modified SBM DEA model that can fully rank all DMUs.
\par

For an evaluated DMU, DEA provides the efficiency (inefficiency) score and a projection point as the improvement target. Both SBM (ERGM) and RM include an efficiency measure, i.e., monotonicity, such that a superior DMU has a better efficiency score than a less effective one, which is a desirable property of the models.   
However, their projection points are not always the closest ones, 
implying a lack of similarity to the evaluated DMU. The concept of similarity of the projection point has been incorporated into prior studies on the least-distance DEA \citep[e.g.,][]{frei1999projections,briec1999holder,ando2012least,aparicio2016survey} and its extensions \citep[e.g.,][]{portela2003finding, ToneSBMMAX,fukuyama2014distance,KAO2022102525}{}. The least-distance projection point is known as the closest target of the evaluated DMU.
Previous studies have attempted to reverse the optimisation of the conventional DEA model to find the closest targets.
For example, \citet{aparicio2007closest} and \citet{ToneSBMMAX} converted the original SBM model into a maximisation model to bring the projection point of the SBM model closer to the evaluated DMU~\citep[see also][]{KAO2021,KAO2022102525}. However, their modified SBM models do not satisfy monotonicity~\citep[see a counterexample of][]{fukuyama2014distance}.
 As argued by~\cite{ando2012least}, not only does the least-distance-based inefficiency measure suffer from the failure of weak monotonicity over the efficient frontier, but so does the additive-form least-distance efficiency measure.
The latter has  an equivalent objective function to  that of the range adjusted measure (RAM) as developed by~\cite{cooper1999ram}.
\par
From a practical perspective, it is necessary to develop a DEA model that can resolve a conflict between the closest targets and strong monotonicity of efficiency measures. 
Strong monotonicity is commonly regarded as a more desirable property than weak monotonicity for a well-defined efficiency measure, particularly in the context of least-distance DEA~\citep[For a detailed discussion, see, e.g.,][]{baek2009relevance,pastor2010relevance}.
In fact, \cite{aparicio2014closest} and \cite{zhu2022determining} incorporated 
 an extended facet approach \citep[e.g.,][]{bessent1988efficiency,green1996efficiency,olesen1996indicators}
into the output-oriented least-distance efficiency measure for determining the closest target from the assessed DMU to the strongly efficient frontier of the extended production possibility set. 
They also demonstrated that strong monotonicity can be ensured under the assumption of the existence of full-dimensional efficient facets~\citep[FDEFs; e.g.,][]{olesen1996indicators}. 
\par
\citet{olesen1996indicators} and \citet{raty2002efficient} have indicated that the extension of FDEFs carries significant implications from an economic perspective. For instance, by extending FDEFs, 
the marginal rates of substitution along the frontier are clearly defined and determined. 
This approach helps to avoid issues associated with a lack of data variation. Consequently, this approach allows for the consideration of a well-defined efficient frontier, 
meaning that we operate not within the standard DEA technology but within an extended production possibility set.
We investigates several technical issues on the extended production possibility set and propose a maximum Russell graph measure (max RM) that provides the closest target for the assessed DMU while satisfying strong monotonicity.
\par
The structure of the paper is as follows:
Section~\ref{S:3} briefly introduces the RM and the least-distance DEA models.
Sections~\ref{S:4} and~\ref{S:5}  develop a max RM with desirable properties.  
Section~\ref{S:6} empirically compares the projection points of a max RM DEA model 
in comparison with other existing non-radial efficiency measures.
Section~\ref{S:7} concludes.

\section{Russell graph measure and least-distance DEA}
\label{S:3}
Denote $P$ as the production possibility set. We assume that $P$ is a polyhedron that includes all convex combinations of the observed input--output vector $(\bm{x}_j,\bm{y}_j)\in \bbR^{m+s}_{++}$ for all $j=1,\dots,n$, while satisfying free disposability~\citep[e.g.,][]{shephard1970theory,fare1985measurement}. 

The efficiency of each DMU is typically measured by comparing it to a dominating projection point on the efficient frontier of the production possibility set. 
Traditional radial efficiency measures developed by~\citet[][]{charnes1978measuring} and~\citet[][]{banker1984some} identify the dominating projection points relative to the weakly efficient frontier of $P$, which is defined as
\begin{equation*}
\partial^w \left(P \right):= \left\{\,
(\bm{x},\bm{y})\in P\mid (\bm{x}',-\bm{y}') < (\bm{x},-\bm{y})\Longrightarrow (\bm{x}',\bm{y}') 
\notin  P\,\right\}. \label{weakfron}
\end{equation*}
Consider a subset of $\partial^w \left(P \right)$ as follows: 
\begin{equation*}
\partial^s \left(P \right):=
\left\{ (\bm{x},\bm{y})\in P \left| 
(\bm{x}',\bm{y}') \not=  (\bm{x},\bm{y}), 
(\bm{x}',-\bm{y}') \leq (\bm{x},-\bm{y})
\Longrightarrow (\bm{x}',\bm{y}')\notin P
\right.
\right\}. \label{parP}
\end{equation*}
The definitions of $\partial^s \left(P\right)$ and $\partial^w \left(P\right)$ imply that $\partial^s \left(P\right) \subseteq \partial^w \left(P\right)$. A DMU that satisfies $(\bm{x}, \bm{y}) \in \partial^w \left(P\right) \setminus \partial^s \left(P\right)$ is said to be weakly efficient, and one that satisfies $(\bm{x}, \bm{y}) \in \partial^s \left(P\right)$ is strongly efficient.
\par
The RM~\citep[][]{fare1985measurement} 
is a non-radial efficiency measure for non-oriented DEA models and the RM DEA model is formulated as follows:
\begin{eqnarray} 
\mbox{$F(\bm{x},\bm{y}):=$}
\min. && \frac{1}{m+s}\left( \sum_{i=1}^m \theta_i +\sum_{r=1}^s \frac{1}{ \phi_r }\right) \label{Q0} \\
\mbox{s.t.} && \left(\sum_{i=1}^m\theta_ix_i\bm{e}_i,\ \sum_{r=1}^s \phi_ry_r \bm{e}_r\right) \in P, \label{Q2} \\
&&0 \le \theta_i  \leq 1\, (i=1,\ldots,m),\ 1 \leq \phi_r \, (r=1,\ldots,s), \label{Q3}
\end{eqnarray}
where $\bm{e}_l$ is the $l^{\mbox{th}}$ unit vector of a proper dimensional space. 
The RM DEA model in equations~\eqref{Q0}--\eqref{Q3} is a nonlinear programming problem and can be solved using a second-order cone programming reformulation~\citep[][]{sueyoshi2007computational}{}{} or a semidefinite programming reformulation~\citep[][]{halicka2018russell}{}{}.   
For an  optimal solution $(\bm{\theta}^*,\bm{\phi}^*)$ of~\eqref{Q0}--\eqref{Q3}, 
the projection point of $(\bm{x},\bm{y})$  is then defined as $\left(\sum_{i=1}^m\theta_i^*x_i\bm{e}_i,\ \sum_{r=1}^s \phi_r^*y_r \bm{e}_r\right)$.
The following theorem ensures that this projection point provided by the RM DEA model lies on the strongly efficient frontier.

\begin{theorem}\label{t2.1}
Suppose that $(\bm{x},\bm{y})\in P\cap\bbR^{m+s}_{++}$ and let $(\bm{\theta}^*,\bm{\phi}^*)$ be an  optimal solution of~\eqref{Q0}--\eqref{Q3}. Then, 
\begin{equation*}
\left(\sum_{i=1}^m\theta_i^*x_i\bm{e}_i,\ \sum_{r=1}^s \phi_r^*y_r \bm{e}_r\right)\in \partial^s(P).\label{eqProj}
\end{equation*}
\end{theorem}
See Theorem 3.1 of  \cite{sekitani2023least} for the proof of Theorem~\ref{t2.1}.
\medskip
\par
For  any $(\bm{x},\bm{y})\in P \cap{\bb R}^{m+s}_{++}$,  $F(\bm{x},\bm{y})$ has a positive value not more than 1. Moreover,  the function $F: P \cap{\bb R}^{m+s}_{++}\to (0,1]$ has the following property: 
\begin{theorem}\label{t2.2}
For each  $(\bm{x},\bm{y})\in P \cap{\bb R}^{m+s}_{++}$,
 \begin{equation}
 F(\bm{x},\bm{y}) >  F(\bar{\bm{x}},\bar{\bm{y}}) \mbox{ if }  (\bm{x},-\bm{y}) \leq  (\bar{\bm{x}},-\bar{\bm{y}}) \mbox{ and }  (\bm{x},-\bm{y}) \not=  (\bar{\bm{x}},-\bar{\bm{y}}). 
 \label{StrMono}
 \end{equation}
\end{theorem}
See Proposition~8 of~ \cite{SUEYOSHI2009764} for the proof of Theorem~\ref{t2.2}.
\medskip
\par
Following~\citet{fare1985measurement}, condition~\eqref{StrMono} of $F$ is a desirable property of an efficiency measure  and 
an efficiency measure satisfying~\eqref{StrMono} is  called strongly monotonic. The efficiency measure $F$ satisfies 
\begin{equation}
 F(\bm{x},\bm{y}) \geq   F(\bar{\bm{x}},\bar{\bm{y}}) \mbox{ if }  (\bm{x},-\bm{y}) \leq  (\bar{\bm{x}},-\bar{\bm{y}}),  
 \label{WeakMono}
 \end{equation} 
 which is weaker than property~\eqref{StrMono}. An efficiency measure satisfying~\eqref{StrMono} is  called weakly monotonic. 
 Therefore, the RM $F$  is strongly and weakly monotonic over $P\cap \bbR^{m+s}_{++}$.   
\par
The projection point achieving the closest target plays an important role in benchmark selection and environment efficiency measures, as shown by~\cite{ruiz2022identifying} and \cite{an2015closest}.  
However, the RM DEA model in~\eqref{Q0}--\eqref{Q3} cannot guarantee the closest target.  
The closest target setting has been widely studied in the least-distance DEA literature, which assumes that the closer projection point is easier to improve for the assessed DMU.
A commensurable inefficiency measure of the least-distance DEA model, which was developed by \citet{aparicio2007closest}, is formulated as 
\begin{eqnarray}
\mbox{$H(\bm{x},\bm{y}):=$}\min. &&  \left\| \left( 1-\theta_1,\cdots,1-\theta_m,\phi_1-1,\cdots,\phi_s-1\right)\right\|_1 \label{l1}\\ 
\mbox{s.t.} && \left(\sum_{i=1}^m\theta_ix_i\bm{e}_i,\ \sum_{r=1}^s \phi_ry_r \bm{e}_r\right) \in \partial^s(P)\label{l2}\\ 
&& 0 \le \theta_i  \leq 1\, (i=1,\ldots,m),\ 1 \leq \phi_r \, (r=1,\ldots,s), \label{l3}
\end{eqnarray}
where $\left\|\cdot\right\|_1$ is the $L_1$ norm, i.e., $\left\|\bm z \right\|_1= \sum_{i=1}^k |z_i|$ for $\bm z = (z_1,\ldots,z_k)$.
Let $(\bm{\theta}',\bm{\phi}')$ be an optimal solution of the model in~\eqref{l1}--\eqref{l3}. The closest target of $(\bm{x},\bm{y})$ is then given by $\left(\sum_{i=1}^m\theta_i'x_i\bm{e}_i,\ \sum_{r=1}^s \phi_r'y_r \bm{e}_r\right)$. 
The projection point $\left(\sum_{i=1}^m\theta_i^*x_i\bm{e}_i,\ \sum_{r=1}^s \phi_r^*y_r \bm{e}_r\right)$ of the RM DEA model in~\eqref{Q0}--\eqref{Q3} is not the closest target to $\partial^s(P)$ if 
\begin{eqnarray*}
    \left\| \left( 1-\theta_1^*,\cdots,1-\theta_m^*,\phi^*_1-1,\cdots,\phi^*_s-1\right)\right\|_1> \\
    \left\| \left( 1-\theta_1',\cdots,1-\theta_m',\phi'_1-1,\cdots,\phi'_s-1\right)\right\|_1, 
\end{eqnarray*}
or, equivalently, $\sum_{r=1}^s \phi_r^*-\sum_{i=1}^m\theta_i^* > \sum_{r=1}^s \phi_r'-\sum_{i=1}^m\theta_i'$.
\par
The least-distance  inefficiency measure $H$ over  $P \cap{\bb R}^{m+s}_{++}$ is represented by the function $H: P \cap{\bb R}^{m+s}_{++}\to [0,\infty)$.
\citet{ando2012least} demonstrated that the inefficiency measure $H$ satisfies neither weak nor strong monotonicity.
\par
To ensure strong monotonicity of inefficiency measures of the least-distance DEA models, 
\cite{aparicio2014closest} adopted an extended facet approach \citep[e.g.,][]{bessent1988efficiency,green1996efficiency,olesen1996indicators} and formulated an output--oriented version
of the model in~\eqref{l1}--\eqref{l3} as 
\begin{eqnarray}
\mbox{$H^+(\bm{x},\bm{y}):=$} \min. & & \left\| \left( \phi_1-1,\cdots,\phi_s-1\right)\right\|_1 \label{ll1}\\
\mbox{s.t.} && \left( \bm{x}, \ \sum_{r=1}^s \phi_ry_r \bm{e}_r\right) \in \partial^s(P_{\scalebox{0.5}{\it EXFA}})\label{ll2}\\ 
 && 1 \leq \phi_r \ (r=1,\ldots,s), \label{ll3}
\end{eqnarray} 
where $P_{\scalebox{0.5}{\it EXFA}}$ represents an extended facet production possibility set (which will be discussed in detail later). 
Then, we have the function $H^+: P \cap{\bb R}^{m+s}_{++}\to [0,\infty)$ as an inefficiency measure of the output-oriented least-distance DEA model.  While $H^+$ satisfies strong monotonicity, one practical disadvantage of the inefficiency measure is that the worst score is preliminarily unknown, whereas any efficiency score of  an efficiency measure is defined  from $0$ to $1$~(see \cite{russell2009axiomatic} for the range of an efficiency measure). As a result, inefficiency measures are not as popular as efficiency measures in DEA applications.

\section{Maximum Russell graph  measure and free-lunch issue in the extended facet production possibility set}
\label{S:4}
To resolve the conflict between a strongly monotonic efficiency measure and  the closest target, we consider a non-oriented DEA model,
\begin{eqnarray}
\max. && \frac{1}{m+s}\left( \sum_{i=1}^m \theta_i +\sum_{r=1}^s \frac{1}{ \phi_r }\right)  \label{p1}\\ 
\mbox{s.t.} && \left(\sum_{i=1}^m \theta_i x_i  \bm{e}_i, \ \sum_{r=1}^s \phi_ry_r \bm{e}_r\right) \in \partial^s(P_{\scalebox{0.5}{\it EXFA}})\label{p2}\\ 
&& 0\leq \theta_i \leq 1 \, (i=1,\ldots,m),\  1\leq \phi_r \ (r=1,\ldots,s),\label{p3}
\end{eqnarray} 
which is a hybrid model of the RM DEA model in~\eqref{Q0}--\eqref{Q3} and 
the least-distance DEA model in~\eqref{l1}--\eqref{l3} with a slight modification by replacing $ \partial^s(P)$ with $ \partial^s(P_{\scalebox{0.5}{\it EXFA}})$. This modification is critical, because without it,
there exists a counterexample in Appendix~\ref{addapendix2} showing that the unmodified efficiency measure fails to satisfy the strong monotonicity. 
\par
As preparation for clarifying the properties of our proposed model, we will now show how to expand   $P_{\scalebox{0.5}{\it EXFA}}$ from $P$ and further discuss some technical details on $P_{\scalebox{0.5}{\it EXFA}}$. 
Assuming variable returns to scale (VRS), we can empirically construct the production possibility set $P$ as follows: $P:=$
\begin{equation}
\left\{
\left(\bm{x},\bm{y}\right) \in {\bb R}^{m+s}_+
\left|\,
\sum_{j=1}^n \lambda_j \bm{x}_j \leq \bm{x}, 
\sum_{j=1}^n \lambda_j \bm{y}_j \geq \bm{y},
\sum_{j=1}^n \lambda_j  = 1, 
\bm{\lambda} \geq \bm{0}
\right.
\right\}.\label{pps}
\end{equation}
A subset  $\cal F$  of $P$ is called a face if there   exists a coefficient vector $(\bm{v},\bm{u},\psi)\in \bbR^{m+s}_+\times \bbR$  such that 
${\cal F}=P \cap \left\{ (\bm{x},\bm{y}) | -\bm{v}^{\top}\bm{x}+\bm{u}^{\top}\bm{y}-\psi=0\right\}$ and 
\[
-\bm{v}^{\top}\bm{x}++\bm{u}^{\top}\bm{y}-\psi\leq 0 \mbox{ for any } (\bm{x},\bm{y})\in P.
\]
If a face $\cal F$ satisfies ${\cal F}\subseteq \partial^s(P)$, then  $\cal F$ is called an efficient face. The efficient face $\cal F$ satisfies 
\[
{\cal F} = P \cap \left\{ (\bm{x},\bm{y}) | -\bm{v}^{\top}\bm{x}+\bm{u}^{\top}\bm{y}-\psi=0\right\}
\]
for some  $(\bm{v},\bm{u})\in \bbR^{m+s}_{++}$  and some $\psi\in \bbR$. See \cite{fukuyama2012efficiency}'s Lemma 2. 
Furthermore, a face $\cal F$ is termed as a FDEF if ${\cal F}\subseteq \partial^s(P)$ and  $\dim {\cal F} = m+ s -1$~\citep[e.g.,][]{olesen1996indicators}. 
\par
Suppose now that $P$ has at least one FDEF. Let ${\cal F}_1,\cdots,{\cal F}_K$ be all the FDEFs of $P$. Then, there exist $(\bm{v}^k,\bm{u}^k)\in \bbR^{m+s}_{++}$ and $\psi^k$ such that 
${\cal F}_k =  P \cap \left\{ (\bm{x},\bm{y}) | -\bm{v}^{k\top}\bm{x}+\bm{u}^{k\top}\bm{y}-\psi^k=0\right\}$ for all $k=1,\ldots,K$. 

In~\citet[][]{olesen1996indicators} and~\citet[][]{aparicio2013well,aparicio2014closest},  
the extended facet approach defines 
an expansion of $P$  based on the FDEFs  as
\begin{eqnarray}
 P_{\scalebox{0.5}{\it EXFA}} :=\left\{ (\bm{x},\bm{y})\in {\bb R}^{m+s}_+  \left|
 -\displaystyle\sum_{i=1}^m v^k_i x_i  + \sum_{r=1}^s u^k_r y_r \leq  \psi^k, k =1,\ldots,K
\right. 
\right\},\label{extPPS16}
\end{eqnarray}
which is called the extended facet production possibility set.

\begin{theorem}\label{thextra}
Assume that all the FDEFs of $P$ defined by \eqref{pps}  are ${\cal F}_1,\cdots,{\cal F}_K$  and let  $J^k=\left\{\, j \, \left| (\bm{x}_j,\bm{y}_j)\in {\cal F}_k\,\right\}\right.$.
The envelopment form of $ P_{\scalebox{0.5}{\it EXFA}}$ defined by \eqref{extPPS16} can be represented as
\begin{eqnarray*}
\left\{ (\bm{x},\bm{y})\in {\bb R}^{m+s}_+  \left|
\sum_{j\in J^k} \lambda_j^k \bm{x}_j  \leq \bm{x},
\sum_{j \in J^k} \lambda_j^k \bm{y}_j \geq \bm{y}, 
\sum_{j\in J^k} \lambda_j^k =1, k =1,\ldots,K
\right\}. \right. 
\end{eqnarray*}
\end{theorem}
See Appendix~\ref{appendixB} for the proof of Theorem~\ref{thextra}. 
 \medskip
\par
The extended facet approach is one of expanding the original production possibility set $P$. \cite{podinovski2004production} shows that the production trade-off assumption expands $P$ and
the expanding set  may violate the no free lunch axiom, which is discussed in \citet{shephard1970theory},~\cite{fare1985measurement}, and~\citet{russell2018technological}. 
Free lunch for a production possibility set  $\cal P$  is characterised as an existence of a certain  input--output vector of $\cal P$ as follows:
\begin{definition}[Free-lunch vector]\label{def1}
If there exists $\bm y\in \bb R_+^s\setminus \{0\}$ such that $(\bm 0,\bm y)\in \cal P$, then the production possibility set $\cal P$ allows free lunch and $(\bm 0,\bm y)$ is called a free-lunch vector. 
\end{definition}
\par
The following example  shows that the set $ P_{\scalebox{0.5}{\it EXFA}} $  allows free lunch. 
The production possibility set $P$ under the VRS assumption has three efficient DMUs: DMU$_A$, DMU$_B$, and DMU$_C$.
\begin{table}[H]
\centering
\tbl{A numeric example}{
\begin{tabular}{crrrrrrr} \hline\hline 
& $A$ & $B$ & $C$ & $D$ & $E$ & $F$ \\ \hline
$x$ &  1   & 2    & 5    & 1   & 1.5 & 20 \\
$y$ &  4   & 5    & 6     & 2     & 2  & 2 \\  \hline\hline 
\end{tabular}}
\label{tab51}
\end{table}
\begin{figure}[H]
\centering
\includegraphics[scale=0.15]{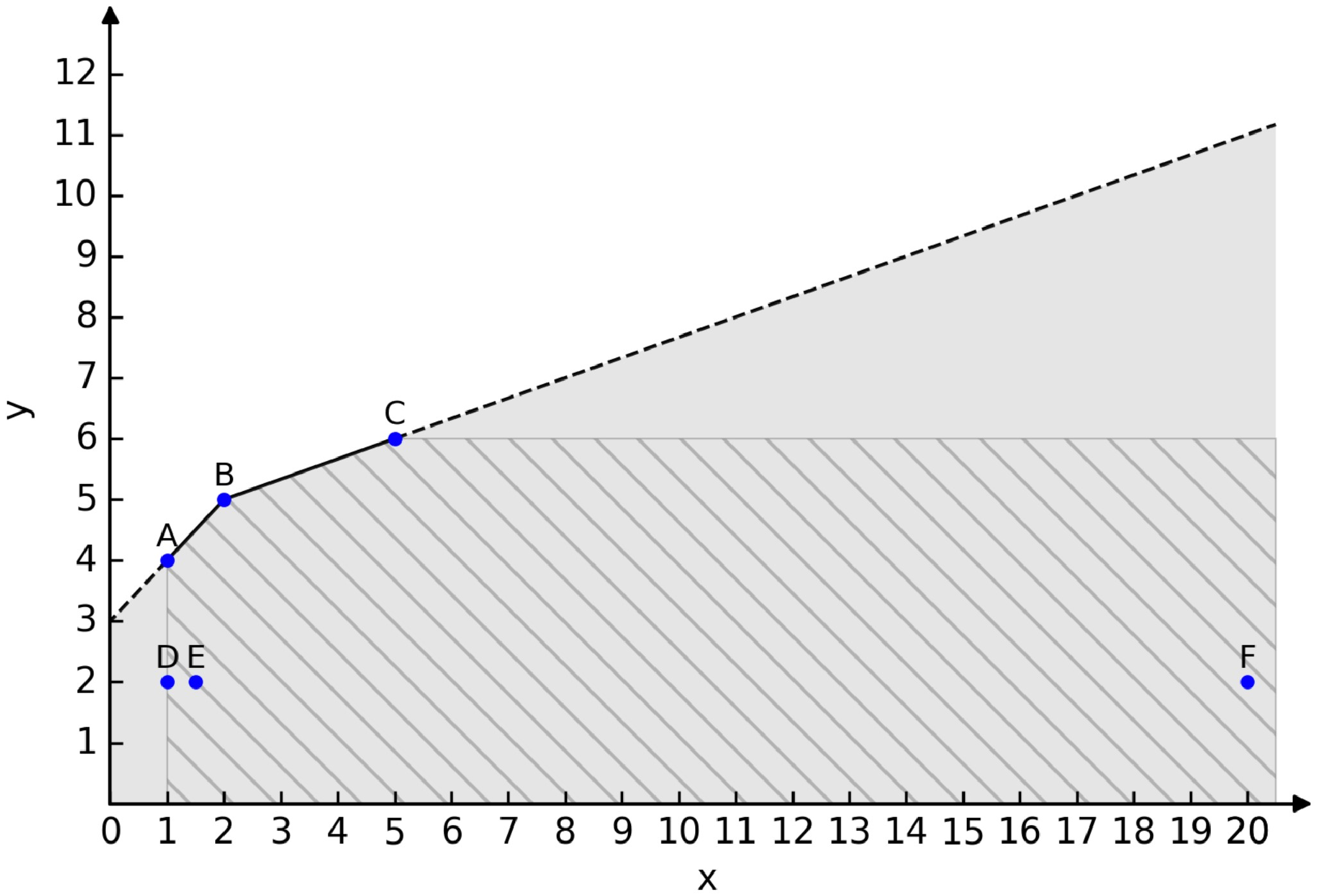}
\caption{The frontier of $P_{\scalebox{0.5}{\it EXFA}}$}\label{fig01}
\end{figure}
As shown in Figure~\ref{fig01},  $P$ is expressed by the gray region filled with shaded lines: 
\begin{align*}
     P=& \left\{ (x,y) \left| 
 \begin{array}{l}   
    \sum_{j=A}^C \lambda_j(x_j,-y_j) \leq (x,-y) , \ \sum_{j=A}^C\lambda_j=1,\\
     y\geq 0,  \     \lambda_j\geq 0 \ (j=A,B,C)
  \end{array}   \right\} \right. 
\end{align*}
and the strongly frontier $\partial^s(P)$ of $P$ consists of two line segments (i.e., the bold line-segment AB and BC): $\partial^s(P)=$
\begin{align*}
   &\left\{(x,y)\mid x-y+3=0,1\le x\le 2\right\}  \cup\left\{(x,y)\mid x-3y+13=0,2\le x\le 5\right\}.
\end{align*}
It follows from $(0,0) \notin P$ that $(0,y) \notin P$ for any $y>0$. Therefore, the set $P$ does not allow free lunch.
\par
The extended facet production possibility set $P_{\scalebox{0.5}{\it EXFA}}$ (i.e., the gray region) can be represented as
\begin{align*}
    P_{\scalebox{0.5}{\it EXFA}}={\Pcon}\cap \bb R_+^{2},
\end{align*}
where $   \Pcon=\left\{(x,y)\left|\  x-y+3\ge0,  \         x-3y+13\geq 0  
    \right.
    \right\}$.
The strongly and weakly efficient frontiers of $P_{\scalebox{0.5}{\it EXFA}}$ are given by
\begin{align*}
\partial^s(P_{\scalebox{0.5}{\it EXFA}})=\partial^s(\Pcon)  \cap  \bbR_+^{2} \mbox{ and }
 \partial^w(P_{\scalebox{0.5}{\it EXFA}})=\partial^s(P_{\scalebox{0.5}{\it EXFA}}) \cup\left\{(0,y)\mid 0\le y\le 3\right\},
\end{align*}
respectively. Since  $(0,y)\in P_{\scalebox{0.5}{\it EXFA}}$ for any positive $y\leq 3$ is a free-lunch vector,  $P_{\scalebox{0.5}{\it EXFA}}$ allows free lunch.
\par

We also demonstrate this issue using a real-world dataset from the target-setting DEA study~\cite{aparicio2007closest}. 
This dataset includes 28 international airlines with two outputs and four inputs; it was originally introduced by~\cite{coelli2002capacity}. 
Under the assumption of VRS, \citet[][]{fukuyama2012decomposing}{}{} verify that  the production possibility set $P$ constructed by the 28 DMUs has one FDEF that is the convex hull 
of six input--output vectors of DMU$_4$, DMU$_8$, DMU$_{13}$, DMU$_{19}$, DMU$_{24}$, and  DMU$_{27}$. 
The supporting hyperplane of $P$  on the six DMUs has the  coefficients
$v^*_1,\ldots,v^*_4$, $u^*_1,u^*_2$ and $\psi^*$ in Table~\ref{Tab3}.
\begin{table}[H]
\centering
\tbl{  The supporting hyperplane of $P$ }{
\begin{tabular}{ccccccc} \hline\hline
$v_1^*$ & $v_2^*$ &  $v_3^*$ & $v_4^*$& $u_1^*$ & $u_2^*$ &{ $\psi^*$ }\\ \hline
0.407493 & 5.69403 & 10.8921 & 1.28423& 0.959148 & 7.22785&{ 15201.7}    \\   \hline\hline 
\end{tabular}}
\label{Tab3}
\end{table}
Since  the FDEF of $P$ is unique,  we have  
\begin{align*}
P_{\scalebox{0.5}{\it EXFA}}
&=\left\{ 
 (\bm{x},\bm{y})\in {\bb R}^{4+2}_+  \left|
-v_1^*x_1-v_2^*x_2-v_3^*x_3-v_3^*x_4
+u_1^*y_1+u_2^*y_2\leq \psi^*
\right\}, \right. 
\end{align*}
which means from $u_1^*+u_2^*<\psi^*$ that  $(0,0,0,0,1,1)\in P_{\scalebox{0.5}{\it EXFA}}$. Hence,  the extended facet production possibility set $P_{\scalebox{0.5}{\it EXFA}}$ allows free lunch.
\par
In general, the extended facet production possibility set $P_{\scalebox{0.5}{\it EXFA}}$ of~\eqref{extPPS16}   allows free lunch if and only if
$\psi^k>0$ for all $k=1,\ldots,K$. 
Hence,  $P_{\scalebox{0.5}{\it EXFA}}$ does not allow free lunch if and only if 
there exists an FDEF having $\psi^k\leq 0$ for some $k \in \{ 1,\ldots,K\}$.

\section{Properties of maximum Russell graph measure}
\label{S:5}
Let $(\bm{\theta},\bm{\phi})$  be a feasible solution of the max RM DEA model in~\eqref{p1}--\eqref{p3}.   
Then, it follows from~\eqref{p3}  that   
 the objective function value of any feasible solution of the DEA model in~\eqref{p1}--\eqref{p3} ranges from 0 to 1.
Let $G(\bm{x},\bm{y})$ be the maximum value of the DEA model in~\eqref{p1}--\eqref{p3}. We have 
\[
0< G(\bm{x},\bm{y}) \leq 1 \mbox{ for all } (\bm{x},\bm{y}) \in   P \cap{\bb R}^{m+s}_{++}.
\]
We therefore call $G(\bm{x},\bm{y})$  a max RM on $P \cap{\bb R}^{m+s}_{++}$. 
\par
The max RM satisfies strong monotonicity. To show this, let  
\begin{eqnarray}
g(\bm{\theta},\bm{\phi}):= \frac{1}{m+s} \left(\sum_{i=1}^m \theta_i+ \sum_{r=1}^s \frac{1}{\phi_r} \right).\label{RM}
\end{eqnarray}  
Define now
 \begin{eqnarray}
G^w(\bm{x},\bm{y}):=\max\left\{\, g(\bm{\theta},\bm{\phi})\,\left|\, 
\begin{array}{c}
 \left(\sum_{i=1}^m \theta_ix_i\bm{e}_i, \sum_{r=1}^s \phi_r y_r \bm{e}_r\right) \in \partial^w(P_{\scalebox{0.5}{\it EXFA}}), \\ 
 \eqref{p3}
\end{array}
 \right\}.\right.\label{maxRMw}
 \end{eqnarray}
Since we have 
  \begin{eqnarray}
G(\bm{x},\bm{y})=\max\left\{\, g(\bm{\theta},\bm{\phi})\,\left|\, 
\begin{array}{c}
 \left(\sum_{i=1}^m \theta_ix_i\bm{e}_i, \sum_{r=1}^s \phi_r y_r \bm{e}_r\right) \in \partial^s(P_{\scalebox{0.5}{\it EXFA}}), \\ 
 ~\eqref{p3}
\end{array}
 \right\},\right.\label{maxRMs}
\end{eqnarray}
and $ \partial^s(P_{\scalebox{0.5}{\it EXFA}})\subseteq  \partial^w(P_{\scalebox{0.5}{\it EXFA}})$, we then obtain the following property:
 \begin{lemma}~\label{Le3.1}
 For any  $(\bm x,\bm y)\in P\cap{\bb R}^{m+s}_{++}$, 
$G^w(\bm{x},\bm{y})\geq G(\bm{x},\bm{y})$. 
 \end{lemma}
See Appendix~\ref{appendixB} for the proof of Lemma~\ref{Le3.1}. 
 \medskip
 \par
 Let $\delta_i:=1-\theta_i$ for all $i=1,\ldots,m$. Then, we can obtain
\[
g(\bm{\theta},\bm{\phi})
= \frac{m}{m+s} + \frac{1}{m+s}  \left( \sum_{i=1}^m (-\delta_i) +\sum_{r=1}^s \frac{1}{\phi_r} \right).
\]
Define $g'(\bm{\delta},\bm{\phi}):= \frac{1}{m+s}  \left( \sum_{i=1}^m (-\delta_i) +\sum_{r=1}^s \frac{1}{\phi_r}\right)$. We, therefore, have  
$g(\bm{\theta},\bm{\phi})= \frac{m}{m+s} +g'(\bm{\delta},\bm{\phi})$ and $g'$ is a quasi-convex function over $\bbR^{m}_{+}\times \bbR^{s}_{++}$ satsifying 
\begin{equation}
g'(\bm{\delta},\bm{\phi}) > g'\left(\bar{\bm{\delta}},\bar{\bm{\phi}}\right) \mbox{ if } (\bm{\delta},-\bm{\phi}) \leq (\bar{\bm{\delta}},-\bar{\bm{\phi}}) \mbox{ and }  (\bm{\delta},-\bm{\phi}) \not= (\bar{\bm{\delta}},-\bar{\bm{\phi}})
\label{decrease}
.\end{equation} 
\textcolor{black}{Let $\bm{e}$ be a vector in a properly dimensional space, where all components are 1.}
The  efficiency measure  $G^w(\bm{x},\bm{y})$ satisfies weak monotonicity \citep{sekitani2023least} and the following statements hold: 
\begin{lemma}~\label{Le3.0}
For any  $(\bm x,\bm y)\in P \cap{\bb R}^{m+s}_{++}$,
\begin{equation}
G^w(\bm{x},\bm{y}) = 
\max\left\{\,
\begin{array}{l}
\max\left\{\left. g(\bm{e}-{(1-\theta'_i)\bm{e}_i},\bm{e}) \right| \, i=1,\ldots.m \right\} \\ 
\max\left\{\left. g(\bm{e},{\bm{e}+(\phi'_r-1)\bm{e}_r}) \right| \, r=1,\ldots.s \right\}
\end{array}
\,\right\}, \label{eqGw}
\end{equation}
where  
\begin{align}
\theta'_i=&\min\left\{ \theta \left|
\begin{array}{l}
\left(\bm{x}-(1-\theta)x_i\bm{e}_i ,\bm{y} \right)\in 
P_{\scalebox{0.5}{\it EXFA}}\\ 
0 \le \theta \leq 1 
\end{array} \right.\right\} \mbox{ for all } i=1,\ldots,m, \label{thetai}\\
\intertext{  and }   
 \phi'_r=&\max\left\{ \phi \left| 
 \begin{array}{l}
\left(\bm{x}, \bm{y}+(\phi-1)y_r\bm{e}_r \right)\in 
P_{\scalebox{0.5}{\it EXFA}}\\ 
1\leq \phi
\end{array}  
 \right.\right\} \mbox{ for all } r=1,\ldots,s.\label{phir}
  \end{align}
 \end{lemma}
  See Appendix~\ref{appendixB} for the proof of Lemma~\ref{Le3.0}. 
 \medskip
 \par
 \begin{lemma}~\label{Le3.2}
 For any $(\bm x,\bm y)\in P\cap{\bb R}^{m+s}_{++}$,  $\bm{\theta}'$ of~\eqref{thetai} and $\bm{\phi}'$ of~\eqref{phir} satisfy
\begin{align*}
\left(\bm{x}-(1-\theta'_i)x_i\bm{e}_i ,\bm{y} \right) \in \partial^s(P_{\scalebox{0.5}{\it EXFA}}) \mbox{ and } \theta'_i > 0 & \mbox{ if } G^w(\bm{x},\bm{y})=  g(\bm{e}-(1-{\theta'_i)\bm{e}_i},\bm{e}) \,\\
\left(\bm{x}, \bm{y}+(\phi'_r-1)y_r\bm{e}_r \right) \in \partial^s(P_{\scalebox{0.5}{\it EXFA}}) &  \mbox{ if } G^w(\bm{x},\bm{y})= g(\bm{e},\bm{e}+({\phi'_r-1)\bm{e}_r}),
\end{align*}
and 
\begin{equation}
G^w(\bm{x},\bm{y})= 
\max\left\{\,
\begin{array}{l}
\max\left\{\, g(\bm{e}-(1-{\theta'_i)\bm{e}_i},\bm{e})\, \left|\, \theta'_i >0 \, \right. \right\}\\ 
\max\left\{\, g(\bm{e},\bm{e}+({\phi'_r-1)\bm{e}_r})\, \left|\, r=1,\ldots,s \,\right.\right\}
\end{array}
\right\} .\label{eqGw1}
\end{equation}
\end{lemma}
See Appendix~\ref{appendixB} for the proof of Lemma~\ref{Le3.2}. 
 \medskip
 \par
Lemma~\ref{Le3.2} implies that $G^w(\bm{x},\bm{y})$ provides the projection point onto $\partial^s(P_{\scalebox{0.5}{\it EXFA}})\setminus \partial^w(P_{\scalebox{0.5}{\it EXFA}})$. 
Therefore, we have $G^w(\bm{x},\bm{y})= G(\bm{x},\bm{y})$ for any  $(\bm{x},\bm{y})\in P\cap{\bb R}^{m+s}_{++}$. Moreover, 
 Theorem~\ref{Th3.3} shows that the efficiency measure  $G(\bm{x},\bm{y})$ satisfies strong monotonicity  on  $P\cap{\bb R}^{m+s}_{++}$.
 \begin{theorem}~\label{Th3.3}
  For any  $(\bm{x},\bm{y})\in P\cap{\bb R}^{m+s}_{++}$,  
$G^w(\bm{x},\bm{y})= G(\bm{x},\bm{y})$. Moreover, the efficiency measure $G(\bm{x},\bm{y})$ is strongly monotonic on $P\cap{\bb R}^{m+s}_{++}$. 
 \end{theorem}
See Appendix~\ref{appendixB} for the proof of Theorem~\ref{Th3.3}. 
 \medskip
 \par
We next show the properties of the projection point provided by the max RM DEA model in~\eqref{p1}--\eqref{p3}. 
The commensurable H\"older input distance function is defined as 
\begin{equation*}
D^-(\bm{x},\bm{y}):=\min\left\{ \sum_{i=1}^m(1- \theta_i) \left|  \left(\sum_{i=1}^m\theta_ix_i\bm{e}_i,\bm{y}\right) \in \partial^w(P_{\scalebox{0.5}{\it EXFA}}),\ \bm{0} \le\bm{\theta} \leq \bm{e}
\right\}\right. 
\end{equation*}
and the commensurable H\"older output distance function is
\begin{equation*}
D^+(\bm{x},\bm{y}):=\min\left\{  \sum_{r=1}^s (\phi_r-1) \left|  \left(\bm{x}, \sum_{r=1}^s \phi_r y_r \bm{e}_r\right) \in \partial^w(P_{\scalebox{0.5}{\it EXFA}}), \bm{e} \leq \bm{\phi} \right\}.\right.
\end{equation*} 
We obtain the following results from Corollary 3 of~\cite{briec1999holder}:  
\begin{align}
&D^-(\bm{x},\bm{y})=1- \max_{i=1,\ldots,m} \min\left\{\, \theta \, \left|\, (\bm{x}-(1-\theta)x_i\bm e_i,\bm{y}) \in P_{\scalebox{0.5}{\it EXFA}} \right\} \right., \label{D-} \\
&D^+(\bm{x},\bm{y})=\min_{r=1,\ldots,s}\max\left\{\, \phi\, \left|\, (\bm{x}, \bm{y}+(\phi-1)y_r\bm{e}_r) \in P_{\scalebox{0.5}{\it EXFA}} \right\} -1 \right.. \label{D+}
\end{align} 
Using $D^-(\bm{x},\bm{y})$ and $D^+(\bm{x},\bm{y})$,  the max RM  can be computed as follows:
 \begin{theorem} \label{Th3.4}
   For any  $(\bm{x},\bm{y})\in P\cap{\bb R}^{m+s}_{++}$, 
\[
G(\bm x,\bm y)= \frac{1}{m+s} \times\max\left\{ 
m+s-D^-(\bm{x},\bm{y}),\ m+s-\frac{D^+(\bm{x},\bm{y})}{1+D^+(\bm{x},\bm{y})}
\right\}.
\]
 \end{theorem}
 See Appendix~\ref{appendixB} for the proof of Theorem~\ref{Th3.4}. 
 \medskip
 \par
A projection point  provided by the commensurable H\"older input distance function  $D^-(\bm{x},\bm{y})$ is called  the   input-oriented  closest target  and 
 a projection point  provided by the commensurable H\"older output distance function  $D^+(\bm{x},\bm{y})$ is called the   output-oriented closest target. 
Theorem~\ref{Th3.4} indicates that a projection point  provided by the max RM DEA model in~\eqref{p1}--\eqref{p3} is either the  input-oriented  closest target or  the output-oriented closest target.
 Furthermore, expression~\eqref{eqGw1} in Lemma~\ref{Le3.2} implies that there exist a positive input element in the projection point. These properties of the projection point provided by the max RM DEA model in~\eqref{p1}--\eqref{p3}
are  summarised as follows:
\begin{corollary}  \label{Co3.5}
Suppose that $(\bm{x},\bm{y})\in P \cap \bbR^{m+s}_{++}$ and let $(\bm{\theta}^*,\bm{\phi}^*)$ be an optimal solution of the max RM DEA model in~\eqref{p1}--\eqref{p3}. Then,
 the projection point $(\sum_{i=1}^m \theta_i^*x_i\bm{e}_i, \sum_{r=1}^s \phi^*_ry_r\bm{e}_r)$ is either the 
  input-oriented closest target or  the output-oriented one. Moreover,  $\sum_{i=1}^m \theta_i^*x_i\bm{e}_i\in \bbR^m_{++}$ and 
   $(\sum_{i=1}^m \theta_i^*x_i\bm{e}_i, \sum_{r=1}^s  \phi^*_ry_r\bm{e}_r)$ 
   is not a free-lunch vector. 
\end{corollary}
See Appendix~\ref{appendixB} for the proof of Corollary~\ref{Co3.5}. 
 \medskip
 \par
 \par
 If $P_{\scalebox{0.5}{\it EXFA}}$ is represented by a Minkowski sum of a convex hull of a
finite number of vertices and a conical hull of a finite number of rays, then two optimisation problems,  
\[
\min\left\{ \theta \left| (\bm{x}-(1-\theta)x_i\bm e_i,\bm{y}) \in P_{\scalebox{0.5}{\it EXFA}} \right\} \right. \mbox{ and }  
\max\left\{ \phi \left| (\bm{x},\bm{y}+(\phi-1)y_r\bm{e}_r) \in P_{\scalebox{0.5}{\it EXFA}} \right\}, \right.
\]
are formulated as linear programming (LP) problems~\citep[see][for details]{sekitani2023least}{}{}, and hence,  $G(\bm{x},\bm{y})$  can be solved by $(m+s)$ LPs.
\section{Comparisons with SBM and its variations}\label{S:6}
\subsection{SBM and maximum SBM}

\citet[][]{pastor1999enhanced}{}{} modified the objective function $g(\bm{\theta},\bm{\phi})$  of the RM DEA model  into $\rho(\bm{\theta},\bm{\phi}):=$
\[
 \frac{\frac{1}{m}\sum_{i=1}^m \theta_i}{ \frac{1}{s}\sum_{r=1}^s  \phi_r } 
\]
and developed an ERGM  DEA model  
$ 
\min \left\{ \left. \rho(\bm{\theta},\bm{\phi}) \right| \,  \eqref{Q2},\ \eqref{Q3}  \right\}, 
$
which is equivalent to the SBM DEA model. 
The following hybrid non-oriented DEA models between the SBM  and the least-distance inefficiency measure are  from~\citet[][]{aparicio2007closest,ToneSBMMAX}{}{} and \cite{aparicio2014closest}:
\begin{align}
    \max&  \left\{\, \rho(\bm{\theta},\bm{\phi})   \left| \left(\sum_{i=1}^m \theta_ix_i\bm{e}_i, \sum_{r=1}^s \phi_r y_r \bm{e}_r\right) \in \partial^s(P), \eqref{Q3} \right\}\right. \label{maxSBM0}\\ 
\intertext{and}
    \max&  \left\{\, \rho(\bm{\theta},\bm{\phi})   \left| \left(\sum_{i=1}^m \theta_ix_i\bm{e}_i, \sum_{r=1}^s \phi_r y_r \bm{e}_r\right) \in \partial^s(P_{\scalebox{0.5}{\it EXFA}}), \eqref{Q3} \right\}\right.,  \label{maxSBM}  
\end{align}
respectively. 
\citet[][]{fukuyama2014distance}{}{} showed that the efficiency measure provided by the hybrid non-oriented DEA model in~\eqref{maxSBM0} violates   weak monotonicity.
\cite{sekitani2023least} later showed that  the efficiency measure provided by the hybrid non-oriented DEA model in~\eqref{maxSBM}  violates strong monotonicity, while \cite{aparicio2014closest} showed that an output-oriented version of~\eqref{maxSBM}   satisfies  strong monotonicity.
An efficiency measure provided by the hybrid non-oriented DEA model in~\eqref{maxSBM} is called a max SBM, which is denoted by $\Gamma(\bm{x},\bm{y})$ for each $(\bm{x},\bm{y})\in P\cap \bbR^{m+s}_{++}$.
\par
For each $(\bm{x},\bm{y})\in P \cap \bbR^{m+s}_{++}$, an optimal solution $(\bm{\theta}^*,\bm{\phi}^*)$ of the max RM DEA model in~\eqref{p1}--\eqref{p3} is 
a feasible solution of the max SBM DEA model in~\eqref{maxSBM}. Since Corollary~\ref{Co3.5} implies $\bm{\theta}^* \in \bbR^{m}_{++}$, we have  
\begin{equation*}
0 < \rho(\bm{\theta}^*,\bm{\phi}^*) \leq \Gamma(\bm{x},\bm{y}) = \rho(\bm{\theta}^{\#},\bm{\phi}^{\#}), \label{eqRHO+}
\end{equation*}
where $(\bm{\theta}^{\#},\bm{\phi}^{\#})$ is an optimal solution of the max SBM DEA model in~\eqref{maxSBM}. 
It follows from $\rho(\bm{\theta}^{\#},\bm{\phi}^{\#})>0$ that  $\sum_{i=1}^m \theta^{\#}_i > 0$ and $\bm{\theta}^{\#}\not= \bm{0}$.
This means that  $\sum_{i=1}^m \theta^{\#}_ix_i\bm{e}_i \not= \bm{0}$ and $\left(\sum_{i=1}^m \theta^{\#}_ix_i\bm{e}_i, \sum_{r=1}^s \phi^{\#}_ry_r\bm{e}_r\right)$ is not a free-lunch vector.
The following corollary implies that the projection point provided by the max SBM DEA model in~\eqref{maxSBM} is not   a free-lunch vector.
\begin{corollary}  \label{Co4.1}
Suppose that $(\bm{x},\bm{y})\in P \cap \bbR^{m+s}_{++}$ and let $(\bm{\theta}^*,\bm{\phi}^*)$ be an optimal solution of the max SBM DEA model in~\eqref{maxSBM}. Then,
  $\sum_{i=1}^m \theta_i^*x_i\bm{e}_i\not= \bm{0}$ and 
   the projection point 
   $\left(\sum_{i=1}^m \theta_i^*x_i\bm{e}_i, \sum_{r=1}^s  \phi^*_ry_r\bm{e}_r\right)$ 
   is not a free-lunch vector. 
\end{corollary}
The projection point provided by the max SBM DEA model in~\eqref{maxSBM} may have an improvement input value of $0$,  while the projection point provided by the max RM DEA model in~\eqref{p1}--\eqref{p3}
avoids a completely vanishing improvement input. 
The existence of a completely vanishing improvement input leads to the violation of  strong monotonicity of the efficiency measure $\Gamma$ as follows: 
\begin{theorem}  \label{Th4.2}
Suppose that $(\bm{x},\bm{y})\in P \cap \bbR^{m+s}_{++}$ and let $(\bm{\theta}^*,\bm{\phi}^*)$ be an optimal solution of the max SBM DEA model in~\eqref{maxSBM}. 
If $\theta_i^*=0$ for some $i \in \{1,\ldots,m\}$, then, 
\[ 
\Gamma(\bm{x},\bm{y}) \leq  \Gamma(\bm{x}+\tau \bm{e}_i,\bm{y})  \mbox{ for any } \tau > 0.
\]
\end{theorem}
See {Appendix}~\ref{appendixB} for the proof of Theorem~\ref{Th4.2}. 
\medskip
\par

Corollaries~\ref{Co4.1} and~\ref{Co3.5} guarantee that both the max SBM and max RM provide projection points that are not a free-lunch vector; however, 
the max SBM  may have a completely vanishing improvement input of the projection point. 
Theorem~\ref{Th4.2} shows that the max SBM $\Gamma$ violates strong monotonicity if a completely vanishing improvement input exists.  
Theorem~\ref{Th3.3} means that the max RM $G$ has a theoretical advantage of strong monotonicity over the max SBM $\Gamma$. 
\subsection{Numerical experiments}
\label{nume}
We now illustrate theoretical results of the max RM DEA model in~\eqref{p1}--\eqref{p2}, the max SBM DEA model in~\eqref{maxSBM}, the SBM DEA model\begin{align}
\min\left\{ \rho(\bm{\theta},\bm{\phi}) \left| \,   \left(\sum_{i=1}^m \theta_ix_i\bm{e}_i, \sum_{r=1}^s \phi_r y_r \bm{e}_r\right) \in P_{\scalebox{0.5}{\it EXFA}}, \eqref{Q3}\right\},\right.
\label{SBM'}
\intertext{and the RM DEA model}
\min\left\{ g(\bm{\theta},\bm{\phi}) \left| \,   \left(\sum_{i=1}^m \theta_ix_i\bm{e}_i, \sum_{r=1}^s \phi_r y_r \bm{e}_r\right) \in P_{\scalebox{0.5}{\it EXFA}}, \eqref{Q3}\right\}\right.
\label{RM'}
\end{align}
using the simple example of Table~\ref{tab51}. We then compare three DEA models~\eqref{p1}--\eqref{p3}, \eqref{maxSBM}, and \eqref{SBM'}
using numerical experiments of four real-world data sets in~\cite{aparicio2007closest},~\citet{wu2011improving},~\citet{juo2016non}, and~\citet{zhu2018simple}. 
The efficiency measures of the SBM~DEA model in~\eqref{SBM'} and the RM DEA model in~\eqref{RM'} are denoted by $\gamma(\bm{x},\bm{y})$ and $f(\bm{x},\bm{y})$ for  each $(\bm{x},\bm{y})\in P\cap \bbR^{m+s}_{++}$, 
repsectively.
\par
Table~\ref{tab2} presents the efficiency scores and projection points of six DMUs in Table~\ref{tab51}. 
The efficiency score of the SBM for DMU$_D$ and DMU$_E$ is the same value of $0$, while DMU$_D$ dominates DMU$_E$.
That is,
\[
0=\gamma(1,2)=\gamma(1.5,2) \mbox{ and } (1,-2) \leq (1.5,-2), (1,-2) \not= (1.5,-2).
\]
This means that the SBM violates strong monotonicity.
The efficiency score of RM for DMU$_D$ and  DMU$_E$ is similar to  the SBM. In fact,   
the efficiency score of RM for DMU$_D$ and  DMU$_E$ is the same value of $1/3$,  while  DMU$_D$ dominates DMU$_E$.
Therefore, the RM violates strong monotonicity. 
The following  three efficiency scores of the max RM DEA model in~\eqref{p1}--\eqref{p3} are consistent with the dominance relationship of the input--output vectors among DMU$_D$, DMU$_E$, and DMU$_F$:
\[
G(1,2) > G(1.5,2) > G(20,2),  
\]
which is also derived from Theorem~\ref{Th3.3}. 

Furthermore, as we have discussed in section~\ref{S:4}, the extended facet production possibility set $P_{\scalebox{0.5}{\it EXFA}}$  
allows free lunch. All  projection points of the SBM and the RM for DMU$_D$ and  DMU$_E$ are a free-lunch vector $(0,3)$. 
For DMU$_D$, DMU$_E$, and DMU$_F$,
both the max SBM and max RM provide projection points that are not a free-lunch vector and which are directly derived from Corollaries~\ref{Co4.1} and~\ref{Co3.5}.

\begin{table}[H]
\tbl{Efficiency scores and projection points}{
\setlength{\tabcolsep}{3.5pt}
\centering
\begin{tabular}{cclcclcclcclcc}
\hline\hline
\multicolumn{2}{c}{\multirow{2}{*}{DMU}} &  & \multicolumn{2}{c}{SBM} &  & \multicolumn{2}{c}{max SBM} &  & \multicolumn{2}{c}{RM} &  & \multicolumn{2}{c}{max RM} \\ \cline{4-5} \cline{7-8} \cline{10-11} \cline{13-14} 
\multicolumn{2}{c}{} &  & $\gamma$ &\small Projection &  & $\Gamma$ &\small Projection    &  & $f$ &\small Projection &  & $G$ &\small Projection    \\ \hline
A                    & (1,4)   &  & 1    & -      &  & 1    & -         &  & 1    & -      &  & 1    & -         \\
B                    & (2,5)   &  & 1    & -      &  & 1    & -         &  & 1    & -      &  & 1    & -         \\
C                    & (5,6)   &  & 1    & -      &  & 1    & -         &  & 1    & -      &  & 1    & -         \\
D                    & (1,2)   &  & 0    & (0,3)  &  & 0.50 & (1,4)     &  & 0.33 & (0,3)  &  & 0.75 & (1,4)     \\
E                    & (1.5,2) &  & 0    & (0,3)  &  & 0.44 & (1.4,4.5) &  & 0.33 & (0,3)  &  & 0.72 & (1.4,4.5) \\
F                    & (20,2)  &  & 0    & (0,3)  &  & 0.18 & (20,11)   &  & 0.25 & (2,5)  &  & 0.59 & (20,11)   \\ \hline\hline
\end{tabular}}
\label{tab2}
\end{table}
\par

In Table~\ref{tab2}, the projection points of all inefficient DMUs are exactly the same in both the max SBM and max RM. 
This is a natural result from Corollaries~\ref{Co4.1} and~\ref{Co3.5} and the fact that there is only one input.
However, do the experiments on a real-world data set, not an artificial one, show consistent projection points for all inefficient DMUs between the max SBM and max RM?
We investigate the consistency using four  real-word  data sets that are available from the existing target setting DEA studies
\citep{aparicio2007closest, juo2016non, zhu2018simple, wu2011improving}. 
The experiment results on the four real-word  data sets are summarised in Tables~\ref{tab51A}--\ref{tab53},  which show the number of DMUs by number of improvement items.

\begin{table}[H]
\tbl{ Projection points of 28 DMUs in \cite{aparicio2007closest} }{
\setlength{\tabcolsep}{7pt}
\centering
\begin{tabular}{lrrrrrr}
\hline\hline
DEA model&\multicolumn{6}{c}{Number of DMUs by number of improvement items}\\ \cline{2-7}
&\multicolumn{1}{c}{0 item}&\multicolumn{1}{c}{1 item}&\multicolumn{1}{c}{2 items}&\multicolumn{1}{c}{3 items}&\multicolumn{1}{c}{4 items}&\multicolumn{1}{c}{5 items}\\ \hline
SBM&6 {\small DMUs}&4 {\small DMUs}&2 {\small DMUs}&7 {\small DMUs}&6 {\small DMUs}&3 {\small DMUs}\\
max SBM&6 {\small DMUs}&20 {\small DMUs}&2 {\small DMUs}&-&-&-\\
max RM&6 {\small DMUs}&22 {\small DMUs}&-&-&-&-\\
\hline\hline
\end{tabular}}
\label{tab51A}
\end{table}

\begin{table}[H]
\tbl{ Projection points of 35 DMUs in \cite{juo2016non} }{
\setlength{\tabcolsep}{7pt}
\centering
\begin{tabular}{lrrrrrr}
\hline\hline
DEA model&\multicolumn{6}{c}{Number of DMUs by number of improvement items}\\ \cline{2-7}
&\multicolumn{1}{c}{0 item}&\multicolumn{1}{c}{1 item}&\multicolumn{1}{c}{2 items}&\multicolumn{1}{c}{3 items}&\multicolumn{1}{c}{4 items}&\multicolumn{1}{c}{5 items}\\ \hline
SBM&7 {\small DMUs}&- &9 {\small DMUs}&8 {\small DMUs}&10 {\small DMUs}&1 {\small DMUs}\\
max SBM&7 {\small DMUs}&26 {\small DMUs}&2 {\small DMUs}&-&-&-\\
max RM&7 {\small DMUs}&28 {\small DMUs}&-&-&-&-\\
\hline\hline
\end{tabular}}
\label{tab510}
\end{table}

\begin{table}[H]
\tbl{ Projection points of 30 { DMUs} in \cite{zhu2018simple} }{
\setlength{\tabcolsep}{11pt}
\centering
\begin{tabular}{lrrrrr}
\hline\hline
DEA model&\multicolumn{5}{c}{Number of DMUs by number of improvement items}\\  \cline{2-6}
&\multicolumn{1}{c}{0 item}&\multicolumn{1}{c}{1 item}&\multicolumn{1}{c}{2 items}&\multicolumn{1}{c}{3 items}&\multicolumn{1}{c}{4 items}\\  \hline
SBM&6 {\small DMUs}&6 {\small DMUs}&2 {\small DMUs}&14 {\small DMUs}&4{\small DMUs}\\ 
max SBM&6 {\small DMUs}&21 {\small DMUs}&3 {\small DMUs}&-&-\\
max RM&6 {\small DMUs}&24 {\small DMUs}&-&-&-\\  \hline \hline
\end{tabular}}
\label{tab52}
\end{table}

\begin{table}[H]
\tbl{ Projection points of 23 DMUs in \cite{wu2011improving} }{
\setlength{\tabcolsep}{8.5pt}
\centering
\begin{tabular}{lrrrrrrrrr}
\hline\hline
DEA model & \multicolumn{9}{c}{Number of DMUs by number of improvement items} \\  \cline{2-10}
& &&\multicolumn{1}{c}{0 item} & & \multicolumn{1}{c}{1 item} & & \multicolumn{1}{c}{2 items} & & \multicolumn{1}{c}{3 items} \\ \hline 
SBM & &&12 {\small DMUs} & & - & & 6 {\small DMUs} & & 5 {\small DMUs} \\
max SBM & &&12 {\small DMUs} & & 10 {\small DMUs} & & 1 {\small DMUs} & & - \\
max RM &&& 12 {\small DMUs} & & 11 {\small DMUs} & & - & & - \\ \hline \hline
\end{tabular}}
\label{tab53}
\end{table}

Theorem~\ref{Th3.4} shows that the number of improvement items for max RM is  either 0 or 1. 
Therefore, DMUs with two or more improvement items in the SBM and max SBM have  projection points that are different from ones provided by the max RM. 
\par
The experiment results in Tables~\ref{tab51A}--\ref{tab53} show that the max SBM and SBM always provide at least one  projection point with two or more improvement items; hence, each experiments on all real-world data sets shows that there exists at least one DMU that has  the inconsistent  projection point between the max SBM and max RM. 
The consistent projection points  of all the DMUs  between the max SBM  and max RM seldom appear in the real world.  
\par

As stated in section~3,  the extended facet production possibility set $P_{\scalebox{0.5}{\it EXFA}}$ generated from the dataset in~\cite{aparicio2007closest} allows free lunch. 
Moreover, $P_{\scalebox{0.5}{\it EXFA}}$ generated from the dataset in~\cite{juo2016non} also allows free lunch. 
In fact, the dataset in~\cite{juo2016non} has three inputs and two outputs, and its extended facet production possibility set is constructed by the supporting hyperplanes of three FDEFs:
\begin{align*}
\scalebox{0.87}{$P_{\scalebox{0.5}{\it EXFA}}=\left\{ 
 (\bm{x},\bm{y})\in {\bb R}^{3+2}_+  \left|
 \begin{array}{ll}
-0.041x_1-6.117x_2-0.003x_3+0.146y_1+0.047y_2\leq 14324.9\\
-0.017x_1-1.748x_2-3.515x_3+0.340y_1+0.011y_2\leq 35073.5\\
-0.025x_1-0.402x_2-3.812x_3+0.391y_1+0.006y_2\leq 40520.3\\
\end{array}
\right\}. \right. $}
\end{align*}
Because all three  supporting hyperplanes have a positive coefficient for the intercept term, 
$P_{\scalebox{0.5}{\it EXFA}}$ allows free lunch. 
\par
The extended facet production possibility sets  generated from the datasets in~\cite{zhu2018simple} and~\citet{wu2011improving} never allow free lunch. 
In fact,  $P_{\scalebox{0.5}{\it EXFA}}$ generated from the dataset in~\cite{zhu2018simple} 
has three FDEFs and at least one of the three FDEFs has   a  {nonpositive coefficient for the intercept term.
The extended facet production possibility set $P_{\scalebox{0.5}{\it EXFA}}$ generated from the dataset in~\cite{,wu2011improving}  
has 18 FDEFs and at least one of the 18 FDEFs has  a nonpositive coefficient for the intercept term.
\par
To further clarify the difference in projection points between the max SBM and max RM, we discuss DMU$_5$ from~\cite{zhu2018simple} as an example. 
This DMU, expressed as $(x_1,x_2,x_3,x_4,y_1,y_2)= (15156,279,1246,2258,12891,$ $599)$, is one of eight in all four real-world datasets with two improvement items. 
The projection point provided by the max SBM is $(15156,279,0,2258,$ $22453.9,599)$, which includes one improvement input $x_3=1246\to0$ and one improvement output $y_1=12891\to22453.9$.
The improvement input $x_3$ is unreasonable due to the vanishing of an input. Moreover, the efficiency score is an evidence of  violating  the strong monotonicity of the corresponding efficiency measure, as per Theorem~\ref{Th4.2}.
However, the projection point provided by the max RM is $(15156,279,1246,2258,36603.5,599)$, which requires only one item to be changed and does not reduce any input to 0.  
The remaining seven DMUs can be explained in the same manner and a completely vanishing of an input is observed in the projections of the max SBM model for all four real-world datasets. Therefore, the theoretical feature of the max RM appears in practice and the max SBM can not become an  alternative target-setting approach of the max RM.

\section{Conclusions}
\label{S:7}
We propose a max RM and its DEA model. Notable properties of the max RM are summarised in Table~\ref{tc}.

\begin{table}[H]
\tbl{Properties of the max RM}{
\setlength{\tabcolsep}{1pt}
\small
\centering
\begin{tabular}{llllll}
\hline\hline
\multirow{2}{*}{DEA } & \multicolumn{1}{c}{Efficiency measure}  &  & \multicolumn{3}{c}{Projection points}                                      \\ \cline{2-2} \cline{4-6} 
\multirow{2}{*}{model }  &
 \multicolumn{1}{c}{\begin{tabular}[c]{@{}c@{}}Strong\\ monotonicity\end{tabular}}& &
  \multicolumn{1}{c}{\begin{tabular}[c]{@{}c@{}}Closest\\ targets\end{tabular}} &
  \multicolumn{1}{c}{\begin{tabular}[c]{@{}c@{}}Maximum \\ improvement\\ number is 1\end{tabular}} &
  \begin{tabular}[c]{@{}c@{}}Not\\ a free-lunch\\ vector\end{tabular} \\ \hline
SBM                        &\quad $\times$ Tables~\ref{tab2}  && \multicolumn{1}{l}{$\times$ Tables~\ref{tab51A}--\ref{tab53}}  & \multicolumn{1}{l}{\quad $\times$ Tables~\ref{tab51A}--\ref{tab53}}  &\quad $\times$  Table~\ref{tab2}\\ 
max SBM                  &\quad $\times$ Thoerem~\ref{Th4.2}  && \multicolumn{1}{l}{$\times$ Tables~\ref{tab51A}--\ref{tab53}}  & \multicolumn{1}{l}{\quad $\times$ Tables~\ref{tab51A}--\ref{tab53}}  &\quad $\circ\ $ Corollary~\ref{Co4.1}  \\ 
max RM                 &\quad $\circ\ $ Theorem~\ref{Th3.3} && \multicolumn{1}{c}{$\circ\ $ Theorem~\ref{Th3.4}} & \multicolumn{1}{c}{\quad $\circ$ Theorem~\ref{Th3.4}} &\quad $\circ\ $ Corollary~\ref{Co3.5} \\ \hline\hline
\end{tabular}}
\label{tc}
\tabnote{$\circ$ indicates that a DEA model satisfies a property. \hfill \quad \hfill\quad \\ 
$\times$ indicates  an opposite case.\hfill \quad \hfill\quad}
\end{table}

All efficiency measures listed in Table~\ref{tc} range from $[0,1]$. Compared with the SBM and max SBM models, the proposed max RM satisfies strong monotonicity and its projection point
of any inefficient DMU is the closest target with exactly one improvement item.
\par
The DEA models under the extended facet production possibility set $ P_{\scalebox{0.5}{\it EXFA}}$  require either enumerating all FDEFs of the original production possibility set  $P$ or a
 mixed integer programming (MIP) solver~\citep{olesen1996indicators}.  Enumerating all the FDEFs has high computational complexity and the MIP solver needs precision parameter tuning  to satisfy the feasibility and  optimality conditions.
 Therefore,  computational complexity and calculation accuracy are serious problems for solving the DEA models under $ P_{\scalebox{0.5}{\it EXFA}}$.
\par

We identified two issues associated with the extended facet production possibility set $P_{\scalebox{0.5}{\it EXFA}}$: the existence of a free lunch vector and the vanishing of an input.
It has been confirmed that both issues occur in four real-world datasets and cannot be ignored practically. The issue of the vanishing of an input occurs independently of the existence of free lunch vectors and is more likely to occur.
The vanishing of an input is a serious concern because it  not only impairs the practicality of projection points but also  undermines the reliability of efficiency measures.
The proposed max RM DEA model can avoid the issue of the vanishing of an input with certainty. It is a future task to verify the existence of DEA models other than the proposed max RM DEA model that can avoid the issue of the vanishing of an input. 

\par
An extended facet production possibility set $ P_{\scalebox{0.5}{\it EXFA}}$ is a kind of an empirical  production possibility set that expands from the original one $P$. 
Imposing the production trade-off assumption~\citep{podinovski2004production} on the production possibility set $P$ expands the original one $P$, which is also an  empirical  expanded  production possibility set   from $P$.
The production possibility under a production trade-off assumption  is often formulated as a dual form of weight restrictions such as assurance region methods. 
The existing DEA models under the  production trade-off assumption associated with assurance region methods never need the enumeration of all FDEFs and MIP solvers~\citep[see,][]{podinovski2007improving}.
This may overcome the problems of computational complexity and calculation accuracy on the empirical expanded  production possibility set. In fact, Theorems~4.1, 4.2, and 4.3 of~\cite{sekitani2023least} show
 that an output-oriented DEA model associated  with  a certain assurance region method does not need to  enumerate all FDEFs  and  solving  $s$ LPs is sufficient to compute the model.
 Therefore, Theorem~\ref{Th3.4} and Corollary~\ref{Co3.5} imply that 
the similar merit of computation could hold for the max RM  DEA model under  the  empirical expanded  production possibility set.

\par
LP solvability  and no free-lunch closest target setting will be beneficial for various DEA applications under more complex structures and circumstances, such as network DEA and imprecise DEA. As part of further research, we aim to modify the max RM DEA model to estimate the best-case efficiency under uncertainty, as discussed in \cite{TOLOO2019313} and \cite{TOLOO2022102583}.

\appendix  

\section{A numerical example of the failure in monotonicity}
\label{addapendix2}
Consider a hybrid non-oriented DEA model of the  RM DEA model in~\eqref{Q0}--\eqref{Q3} and the least-distance DEA model in~\eqref{l1}--\eqref{l3} as follows:
\begin{eqnarray}\label{am01}
\mbox{ $M(\bm{x},\bm{y}):=$}\max. && \frac{1}{m+s}\left( \sum_{i=1}^m \theta_i +\sum_{r=1}^s \frac{1}{ \phi_r }\right) \\ 
\mbox{s.t.} && \left(\sum_{i=1}^m \theta_i x_i  \bm{e}_i, \ \sum_{r=1}^s \phi_ry_r \bm{e}_r\right) \in \partial^s(P)\\\label{am2}
&& 0\leq \theta_i \leq 1 \, (i=1,\ldots,m),\  1\leq \phi_r \ (r=1,\ldots,s).
\label{am3}
\end{eqnarray} 
 A numerical example of Table~\ref{tab01b} shows that the DEA model in~\eqref{am01}--\eqref{am3}  fails to satisfy strong monotonicity
 of the efficiency measure on $P \cap \bbR^{m+s}_{++}$.

\begin{table}[H]
\centering
\tbl{A counterexample on $P$}{
\begin{tabular}{cccccccccccccc}
\hline\hline
\multirow{3}{*}{DMU} &
  \multirow{3}{*}{$x_1$} &
  \multirow{3}{*}{$x_2$} &
  \multirow{3}{*}{$y_1$} &
  \multirow{3}{*}{Status} &
  \multicolumn{4}{c}{RM} &
   &
  \multicolumn{4}{c}{DEA model in~\eqref{am01}--\eqref{am3}} \\ \cline{6-9} \cline{11-14} 
    &    &    &    &        & Eff.    & \multicolumn{3}{c}{Target} &  & Eff.    & \multicolumn{3}{c}{Target} \\
    &    &    &    &        & Score   & $x_1$   & $x_2$   & $y_1$  &  & Score   & $x_1$   & $x_2$   & $y_1$  \\ \hline
$A$ & 1  & 1  & 8  & Eff.   & 1       & -       & -       & -      &  & 1       & -       & -       & -      \\
$B$ & 10 & 5  & 10 & Eff.   & 1       & -       & -       & -      &  & 1       & -       & -       & -      \\
$C$ & 5  & 10 & 10 & Eff.   & 1       & -       & -       & -      &  & 1       & -       & -       & -      \\
$D$ & 10 & 1  & 8  & Ineff. & $7/10$  & 1       & 1       & 8      &  & $7/10$  & 1       & 1       & 8      \\
$E$ & 10 & 5  & 8  & Ineff. & $13/30$ & 1       & 1       & 8      &  & $14/15$ & 10      & 5       & 10     \\ \hline\hline
\end{tabular}}
\label{tab01b}
\end{table}
In Table~\ref{tab01b}, we calculate the efficiency scores of the RM and DEA model in~\eqref{am01}--\eqref{am3} under the assumption of VRS for five DMUs (A, B, C, D, and E). 
 Three DMUs, DMU$_A$, DMU$_B$, and DMU$_C$ are efficient. 
Both DMU$_D$ and DMU$_E$ are inefficient, with DMU$_D$ performing better than DMU$_E$ because DMU$_D$ uses less input amount $x_2$ than DMU$_E$ does to produce the same level of output. 
However, the DEA model in~\eqref{am01}--\eqref{am3} results in a smaller efficiency score for DMU$_D$ than for DMU$_E$ (i.e., $\frac{7}{10}<\frac{14}{15}$), indicating that the DEA model in~\eqref{am01}--\eqref{am3}  fails to satisfy strong monotonicity
 of efficiency measures over $P\cap \bbR^{m+s}_{++}$.
\par
The above observation can be explained as follows:
Since
\begin{align*}
    (x_{1A}, x_{2A}, -y_{1A})=(1,1,-8)&\leq  (10,1,-8)= (x_{1D}, x_{2D}, -y_{1D})
    \intertext{and}
    (x_{1B}, x_{2B}, -y_{1B}) = (10,5,-10)&\leq (10,5,-8)= (x_{1E}, x_{2E}, -y_{1E}),
\end{align*}
DMU$_D$ and DMU$_E$ are inefficient. There exists a supporting hyperplane $2x_1+2x_2-13y_1=100$ of $P$ passing on  DMU$_A$, DMU$_B$, and DMU$_C$. Since the coefficient vector of the supporting hyperplane,
$(2,2,13)$,    is positive,
DMU$_A$, DMU$_B$, and DMU$_C$ are efficient. Three input--output vectors of DMU$_A$, DMU$_B$, and DMU$_C$ are affinely independent,  and hence, 
the convex hull of the three DMUs are  the FDEF of $P$.  Consequently, the strongly efficient frontier $\partial^s(P)$ is represented by the following convex hull:
\begin{align*}
    \left\{(x_1,x_2,y_1)\left|
    \begin{array}{ll}
      (x_1,x_2,y_1)=(1,1,8)\lambda_A+(10,5,10)\lambda_B+(5,10,10)\lambda_C \\
      \lambda_A+\lambda_B+\lambda_C=1,\      \lambda_A,\lambda_B,\lambda_C\ge 0
    \end{array}\right\}.\right.
\end{align*}
The frontier of the corresponding $P$ is shown in Figure~\ref{fig02}, where the FDEF is highlighted in orange.
\begin{figure}[H]
\centering
\includegraphics[scale=0.42]{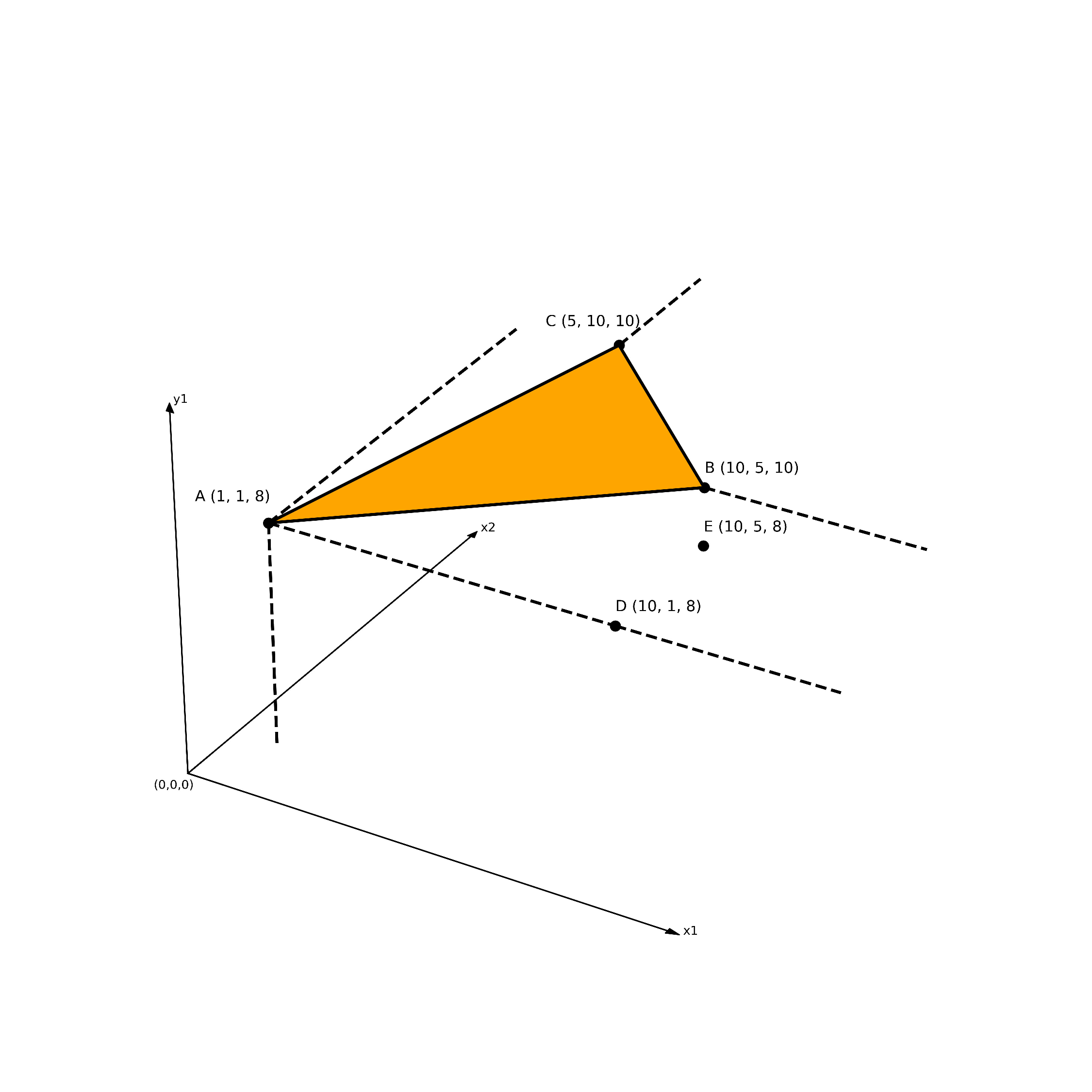}
\caption{The frontier of $P$ having the FDEF  represented in orange}\label{fig02}
\end{figure}
It follows from  $(x_{1A}, x_{2A}, y_{1A})\in \partial^s(P)$ that $M(x_{1A}, x_{2A}, y_{1A})=1$. In the same manner,  we have $M(x_{1B}, x_{2B}, y_{1B})=M(x_{1C}, x_{2C}, y_{1C})=1.$
From  Figure~\ref{fig02}, we see that DMU$_A$ is a unique efficient input--output vector dominating DMU$_D$. 
In fact,   the model in~\eqref{am01}--\eqref{am3} for DMU$_D$ can be
represented by
\begin{align}\label{zzz}
\max. &\  \frac{1}{3}\left( \frac{\lambda_A+10\lambda_B+5\lambda_C}{10}+\frac{\lambda_A+5\lambda_B+10\lambda_C}{1}+\frac{1}{\frac{8\lambda_A+10\lambda_B+10\lambda_C}{8}}\right) \\ 
\mbox{s.t.} &\  
\lambda_A+10\lambda_B+5\lambda_C \le 10 \label{zzA}\\
&\ \lambda_A+5\lambda_B+10\lambda_C\le 1 \label{zzB}\\
&\ 8\lambda_A+10\lambda_B+10\lambda_C\ge 8 \label{zzC}\\
&\ \lambda_A+\lambda_B+\lambda_C=1 \label{zzOne}\\
&\ \lambda_A,\lambda_B,\lambda_C\ge 0.\label{zzD}
\end{align} 
It follows from \eqref{zzB}, \eqref{zzOne}, and \eqref{zzD} that the model in \eqref{zzz}--\eqref{zzD} for  DMU$_D$ has a unique  feasible solution, $(\lambda_A,\lambda_B,\lambda_C)=(1,0,0)$.
Therefore, the optimal solution of  the model in \eqref{zzz}--\eqref{zzD} for  DMU$_D$ is also $(\lambda_A,\lambda_B,\lambda_C)=(1,0,0)$ and the optimal value,  $M(x_{1D},x_{2D},y_{1D})$, is
\[
\frac{1}{3}\left( \frac{1+0+0}{10}+\frac{1+0+0}{1}+\frac{1}{\frac{8+0+0}{8}}\right)=\frac{1}{3}\left( \frac{1}{10} + 1 + 1\right)=\frac{1}{3}\cdot\frac{21}{10}=\frac{7}{10}.
\]

Subsequently, we estimate the lower bound of the 
efficiency level for DMU$_E$. 
The model in~\eqref{am01}--\eqref{am3} for DMU$_E$ can be
represented by
\begin{align}\label{zzzz}
\max. &\  \frac{1}{3}\left( \frac{\lambda_A+10\lambda_B+5\lambda_C}{10}+\frac{\lambda_A+5\lambda_B+10\lambda_C}{5}+\frac{1}{\frac{8\lambda_A+10\lambda_B+10\lambda_C}{8}}\right) \\ 
\mbox{s.t.} &\  
\lambda_A+10\lambda_B+5\lambda_C \le 10 \label{zzzA}\\
&\ \lambda_A+5\lambda_B+10\lambda_C\le 5 \label{zzzB}\\
&\ 8\lambda_A+10\lambda_B+10\lambda_C\ge 8 \label{zzzC}\\
&\ \lambda_A+\lambda_B+\lambda_C=1 \label{zzzOne}\\
&\ \lambda_A,\lambda_B,\lambda_C\ge 0.\label{zzzD}
\end{align} 
Because $(\lambda_A,\lambda_B,\lambda_C)=(0,1,0)$ is a feasible solution of~\eqref{zzzz}--\eqref{zzzD}, we have
\begin{align*}
    M(x_{1E},x_{2E},y_{1E})\ge \frac{1}{3}\left( 1+1+\frac{8}{10}\right)=\frac{14}{15}. 
\end{align*}
Since $M(x_{1D},x_{2D},y_{1D})=\frac{14}{20}<\frac{14}{15}\leq M(x_{1E},x_{2E},y_{1E})$, the efficiency measure $M$ fails to satisfy strong monotonicity.
\par
Finally, we will show $M(x_{1E},x_{2E},y_{1E})=\frac{14}{15}.$ It follows from $\lambda_B=1-(\lambda_A+\lambda_C)$ that the model  in \eqref{zzzz}--\eqref{zzzD} for DMU$_E$ is 
\begin{align}\label{xx}
\max. &\ \frac{1}{3} \left(2- \frac{17}{10} \lambda_A +\frac{1}{2} \lambda_C + \frac{8}{10-2\lambda_A}\right) \\ 
\mbox{s.t.} &\  
-9\lambda_A-5\lambda_C \le 0 \label{xxA}\\
&\ -4\lambda_A+5\lambda_C\le 0 \label{xxB}\\
&\ \lambda_A \le 1 \label{xxC}\\
&\ 0 \leq  \lambda_A + \lambda_C\leq 1 \label{xxOne}\\
&\ \lambda_A\geq 0\ ,\lambda_C\ge 0.\label{xxD}
\end{align} 
All constraints of   \eqref{xxA}--\eqref{xxD}  are equivalently reduced to 
\begin{align}\label{xxx}
0\leq \lambda_C \leq \frac{4}{9} \mbox{ and } \frac{5}{4}\lambda_C \leq \lambda_A \leq 1-\lambda_C. 
\end{align} 
Therefore, the model in \eqref{zzzz}--\eqref{zzzD} for DMU$_E$ is 
\begin{align}
\max\left\{  \left. 
\frac{1}{3}\left( 2- \frac{17}{10} \lambda_A +\frac{1}{2} \lambda_C + \frac{8}{10-2\lambda_A} \right) \right|\ 
\begin{array}{l}
0\leq \lambda_C \leq \frac{4}{9} \\
 \frac{5}{4}\lambda_C \leq \lambda_A \leq 1-\lambda_C
 \end{array}
  \right\}. \label{maxSimple}
\end{align}
For  four positive numbers $a$, $b$, $c$ and $\beta$,   and a nonnegative number  $\alpha < \beta $, let  
\[
\phi(t) := -a t + \frac{b}{c-t} \mbox{ for any } t\in [\alpha,\beta]
\]
and assume $\frac{b}{a}<\left( c-1\right)^2$ and $c>1\geq \beta$. Then, 
for any $t \leq \beta$, the derivative of $\phi$ satisfies
\begin{align*}
{\phi'(t)} = b\left( -\frac{a}{b}+ \frac{1}{(c-t)^2}\right)  \leq b \left( -\frac{a}{b}+ \frac{1}{(c-\beta)^2}\right) \leq  b \left( -\frac{a}{b}+ \frac{1}{(c-1)^2}\right)<0
\end{align*}
and $\phi$ is a decreasing function on $ (-\infty,\beta]$. Therefore, 
\begin{align}
\max_{\alpha \leq t \leq \beta} \phi(t) = \phi(\alpha). \label{maxAlpha}
\end{align}
It follows from \eqref{maxSimple} and  \eqref{maxAlpha} that 
\begin{align*} 
3\times M(x_{1E},x_{2E},y_{1E})
&=\max_{0\leq \lambda_C \leq \frac{4}{9}} 
\left\{  
\left. 
2+ \frac{1}{2} \lambda_C + \max\
\left\{ 
-\frac{17}{10} \lambda_A + \frac{4}{5-\lambda_A} 
\right|  
\frac{5}{4}\lambda_C \leq \lambda_A \leq 1-\lambda_C
\right\}\right\}\\
&=\max_{0\leq \lambda_C \leq \frac{4}{9}} 
\left\{  
2+ \frac{1}{2} \lambda_C  - \frac{17}{8}\lambda_C + \frac{16}{20-5\lambda_C}
\right\}\\
&= \max_{0\leq \lambda_C \leq \frac{4}{9}} 
\left\{  
2   - \frac{13}{8}\lambda_C + \frac{\frac{16}{5}}{4-\lambda_C} 
\right\}=\left( 2 +  \frac{16}{20} \right) = \frac{14}{5}=3\times \frac{14}{15}.
\end{align*} 
This means $M(x_{1E},x_{2E},y_{1E})=\frac{14}{15}.$

\section{Proofs}\label{appendixB}
\label{proofs}

\noindent
\textbf{Theorem~\ref{thextra}}

\begin{proof}
 Let 
 \begin{align*}
{\cal  H}^k_-&:=\left\{\, (\bm{x},\bm{y})\, \left|\,
 -\bm{v}^{k\top}\bm{x}+\bm{u}^{k\top}\bm{y}\leq \psi^k
\, \right\}\right.
\intertext{and}
{\cal F}^k_-&:=\left\{\, (\bm{x},\bm{y})\, \left|\,
\sum_{j\in J^k} \lambda_j \bm{x}_j \leq \bm{x},\  \sum_{j\in J^k} \lambda_j\bm{y}_j \geq \bm{y}, \sum_{j \in J^k} \lambda_j=1
\, \right\}\right.
 \end{align*}
 for all $k=1,\ldots,K$, where  $(\bm{v}^k,\bm{u}^k)\in \bbR^{m+s}_{++}$ and $\psi^k$ satisfy
${\cal F}_k =  P \cap \left\{ (\bm{x},\bm{y}) | -\bm{v}^{k\top}\bm{x}+\bm{u}^{k\top}\bm{y}-\psi^k=0\right\}$.
Moreover, without loss of generality we assume $\sum_{i=1}^m v_i^k + \sum_{r=1}^s u_r^k=1$. Then,  
it suffices to prove ${\cal  H}^k_-={\cal F}^k_-$ for all $k=1,\ldots,K$.
\par
Firstly, we will show  ${\cal  H}^k_-\subseteq {\cal F}^k_-$. Consider a linear programing problem
\begin{align}
&\max\left\{ -\bm{v}^{\top}\bm{x}+\bm{u}^{\top}\bm{y}-\psi\left| 
\begin{array}{l}
-\bm{v}^{\top}\bm{x}_j+\bm{u}^{\top}\bm{y}_j-\psi=0,\  \forall j\in J^k \\
\sum_{i=1}^m v_i + \sum_{r=1}^s u_r=1, \bm{v}\geq \bm{0},\bm{u}\geq \bm{0}\end{array}
\right\}\right.\label{PrimalLP}
\intertext{ and its dual problem}
&\min\left\{ \theta \left| 
\begin{array}{l}
\sum_{j \in J^k} \lambda_j \bm{x}_j -\theta \bm{e} \leq \bm{x}, 
\sum_{j \in J^k} \lambda_j \bm{y}_j + \theta\bm{e}\geq \bm{y}, \sum_{j \in J^k } \lambda_j=1
\end{array}
\right\}.\right.\label{DualLP}
\end{align}
The primal problem~\eqref{PrimalLP} has the unique optimal solution  $(\bm{v}^k,\bm{u}^k,\psi^k)$ since ${\cal F}^k$ is an FDEF.
Suppose that $(\bm{x},\bm{y})\in {\cal  H}^k_-$. Then, the optimal value of~\eqref{PrimalLP} is  $-\bm{v}^{k\top}\bm{x}+\bm{u}^{k\top}\bm{y}-\psi^k\leq 0$.
 Let $(\bm{\theta}^*,\bm{\lambda}^*)$ be an optimal solution of \eqref{DualLP}. Then, we have the optimal value $\theta^*\leq 0$ and 
\begin{align*}
&\bm{x}\geq \sum_{j \in J^k} \lambda_j^* \bm{x}_j -\theta^* \bm{e} \geq  \sum_{j \in J^k} \lambda_j^* \bm{x}_j, \\
&\bm{y} \leq \sum_{j \in J^k} \lambda_j^* \bm{y}_j + \theta^*\bm{e}\leq  \sum_{j \in J^k} \lambda_j^* \bm{y}_j  \\
&\sum_{j \in J^k } \lambda_j^*=1.
\end{align*}
Therefore, $(\bm{x},\bm{y})\in {\cal  F}^k_-$.
\par
Secondly, we will show ${\cal  F}^k_-\subseteq {\cal H}^k_-$. Suppose  $(\bm{x},\bm{y})\in {\cal  F}^k_-$. Then, there exists $\bm{\lambda}$ such that $\sum_{j \in J^k} \lambda_j \bm{x}_j \leq \bm{x}, 
\sum_{j \in J^k} \lambda_j \bm{y}_j \geq \bm{y}$ and $\sum_{j \in J^k } \lambda_j=1. $  Since $-\bm{v}^{k\top}\bm{x}_j+\bm{u}^{k\top}\bm{y}_j-\psi^k=0$ for all $j  \in J^k$ and $(\bm{v}^k,\bm{u}^k)\in \bbR^{m+s}_{++}$, it follows that 
\begin{align*}
-\bm{v}^{k\top}\bm{x}+\bm{u}^{k\top}\bm{y}-\psi^k &\leq -\bm{v}^{k\top}\left( \sum_{j \in J^k} \lambda_j\bm{x}_j \right) +\bm{u}^{k\top}\left( \sum_{j \in J^k}\lambda_j \bm{y}_j \right) -\psi^k\left( \sum_{j \in J^k}\lambda_j  \right)\\
                                                                  & =  - \sum_{j \in J^k} \lambda_j   \left( \bm{v}^{k\top} \bm{x}_j  +\bm{u}^{k\top} \bm{y}_j -\psi^k\right) =0.
\end{align*}
Therefore, we have  $(\bm{x},\bm{y})\in {\cal  H}^k_-$.
\par

 \end{proof}

\noindent
\textbf{Lemma~\ref{Le3.1}}
 \begin{proof}\label{proofLe3.1}
 It follows from $\partial^s(P_{\scalebox{0.5}{\it EXFA}}) \subseteq \partial^w(P_{\scalebox{0.5}{\it EXFA}})$ that 
 an optimal solution of the max RM DEA model in~\eqref{p1}--\eqref{p3} is a feasible solution of the DEA model in~\eqref{maxRMw}. 
 Therefore, the optimal value  of the DEA model in~\eqref{maxRMw} is not less than that of the max RM DEA model in~\eqref{p1}--\eqref{p3}, and hence,  
  $G^w(\bm{x},\bm{y})\geq G(\bm{x},\bm{y})$.
  \end{proof}

\noindent
\textbf{Lemma~\ref{Le3.0}}
\begin{proof}\label{proofLe3.0}
Choose any $(\bm x,\bm y)\in P\cap{\bb R}^{m+s}_{++}$ and let 
$\delta'^{-}_i:= 1-\theta'_i$ for all $i=1,\ldots,m$ and $\delta'^{+}_r:=\phi'_r-1$ for all $r=1,\ldots,s$.
Let 
\begin{equation}
\Delta:= \left\{\, 
(\bm{\delta}^-,\bm{\delta}^+)\,
\left| \,
\begin{array}{l}
 (\bm{\delta}^-,\bm{\delta}^+) = \sum_{i=1}^m \lambda_i(\delta'^{-}_i\bm{e}_i,\bm{0})+
\sum_{r=1}^s \lambda_{m+r} (\bm{0},\delta'^{+}_r\bm{e}_r)\\
\sum_{j=1}^{m+s} \lambda_j=1, \ \bm{\lambda} \geq \bm{0}
\end{array}
\right\}\right.
\end{equation}
and
\begin{equation}
\Omega:= \left.\left\{\, 
\left(\bm{x}-\sum_{i=1}^m \delta^{-}_ix_i\bm{e}_i, \bm{y}+\sum_{r=1}^{s}\delta_r^{+}y_r\bm{e}_r\right)\,\right|\,
(\bm{\delta}^-,\bm{\delta}^+)\in \Delta 
\right\}.
\end{equation}
Since $(\bm{x}-\delta'^{-}_ix_i\bm{e}_i,\bm{y})\in \Pexfa$ for all $i=1,\ldots,m$ and
$(\bm{x},\bm{y}+\delta'^{+}_ry_r\bm{e}_r)\in \Pexfa$ for all $r=1,\ldots,s$,
we have $\Omega \subseteq \Pexfa$.
\par
Let $(\hat{\bm{\theta}},\hat{\bm{\phi}})$  be a feasible solution of the model in~\eqref{maxRMw} and let 
$\hat{\delta}^{-}_i:= 1-\hat{\theta}_i$ for all $i=1,\ldots,m$ and $\hat{\delta}^{+}_r:=\hat{\phi}_r-1$ for all $r=1,\ldots,s$.. 
Then, $(\hat{\bm{\delta}}^{-},\hat{\bm{\delta}}^+)\in \bbR^{m+s}_{+}$ satisfies
\begin{equation}
\left(\bm{x}-\sum_{i=1}^m \hat{\delta}^{-}_ix_i\bm{e}_i, \bm{y}+\sum_{r=1}^{s}\hat{\delta}_r^{+}y_r\bm{e}_r\right)\in \partial^w\left( \Pexfa \right).
\end{equation}
Since any $(\bm{x}'',\bm{y}'')\in  \partial^w\left( \Pexfa \right)$ with 
$(\bm{x}'',-\bm{y}'')\leq  (\bm{x},-\bm{y})$ satisfies 
$\alpha (\bm{x}''-\bm{x},\bm{y}''-\bm{y})+ (\bm{x},\bm{y})\in \Omega$ for some non-negative number $\alpha \leq 1$,
there exists a non-negative $\hat{\alpha} \leq 1$ such that $
 \left(\bm{x}-\hat{\alpha}\sum_{i=1}^m \hat{\delta}^{-}_ix_i\bm{e}_i, \bm{y}+
\hat{\alpha}\sum_{r=1}^{s}\hat{\delta}_r^{+}y_r\bm{e}_r\right) \in \Omega$, and hence,
$\hat{\alpha}(\hat{\bm{\delta}}^{-},\hat{\bm{\delta}}^+)\in \Delta$.
\par
It follows from the definition in~\eqref{RM} of $g$ that the function $g'$ is 
quasi-convex over $\bbR^{m+s}_{+}$ satisfying~\eqref{decrease}, and hence, 
 \begin{align}
 g'\left(\hat{\alpha}\hat{\bm{\delta}}^-,\bm{e}+\hat{\alpha}\hat{\bm{\delta}}^+\right)&\leq \max\left\{  
g'\left(\bm{\delta}^-,\bm{\delta}^+\right) \left| 
 \left(\bm{\delta}^-,\bm{\delta}^+\right)\in \Delta \right\}\right.  \label{eqAA1} \\
& =
\max\left\{
\begin{array}{l}
\max\left\{\left. g'(\delta'^{-}\bm{e}_i,\bm{e}) \right| \, i=1,\ldots.m \right\} \\ 
\max\left\{\left. g'(\bm{0},\bm{e}+\delta'^{+}_r\bm{e}_r) \right| \, r=1,\ldots.s \right\}
\end{array}
\right\}. \label{eqAA2}
\end{align}
Therefore, for any  feasible solution  $(\hat{\bm{\theta}},\hat{\bm{\phi}})$ of the model in~\eqref{maxRMw}, 
it follows from~\eqref{eqAA1} and \eqref{eqAA2} that 
\color{black}
\begin{eqnarray*}
g\left(\hat{\bm{\theta}},\hat{\bm{\phi}}\right)&=&
\frac{m}{m+s}+ g'\left(\hat{\bm{\delta}}^-,\bm{e}+\hat{\bm{\delta}}^+\right)\leq \frac{m}{m+s}+
g'\left(\hat{\alpha}\hat{\bm{\delta}}^-,\bm{e}+\hat{\alpha}\hat{\bm{\delta}}^+\right)\\
&\leq&
{\frac{m}{m+s}+\max\left\{  
g'\left(\bm{\delta}^-,\bm{\delta}^+\right) \left| 
 \left(\bm{\delta}^-,\bm{\delta}^+\right)\in \Delta \right\}\right.}\\
&=&
\frac{m}{m+s}+
\max\left\{
\begin{array}{l}
\max\left\{\left. g'(\delta'^{-}\bm{e}_i,\bm{e}) \right| \, i=1,\ldots.m \right\} \\ 
\max\left\{\left. g'(\bm{0},\bm{e}+\delta'^{+}_r\bm{e}_r) \right| \, r=1,\ldots.s \right\}
\end{array}
\right\}\\
&=& 
\max\left\{
\begin{array}{l}
\max\left\{\left. g(\bm{e}-\delta'^{-}\bm{e}_i,\bm{e}) \right| \, i=1,\ldots.m \right\} \\ 
\max\left\{\left. g(\bm{e},\bm{e}+\delta'^{+}_r\bm{e}_r) \right| \, r=1,\ldots.s \right\}
\end{array}
\right\}\\
&=& 
\max\left\{
\begin{array}{l}
\max\left\{\left. g(\bm{e}-(1-\theta'_i)\bm{e}_i,\bm{e}) \right| \, i=1,\ldots.m \right\} \\ 
\max\left\{\left. g(\bm{e},\bm{e}+(\phi'_r-1)\bm{e}_r) \right| \, r=1,\ldots.s \right\}
\end{array}
\right\}.
\end{eqnarray*}
This means
\[
G^w(\bm{x},\bm{y})=\max\left\{
\begin{array}{l}
\max\left\{\left. g(\bm{e}-(1-\theta'_i)\bm{e}_i,\bm{e}) \right| \, i=1,\ldots.m \right\} \\ 
\max\left\{\left. g(\bm{e},\bm{e}+(\phi'_r-1)\bm{e}_r) \right| \, r=1,\ldots.s \right\}
\end{array}
\right\}
\] 
since $(\bm{e}-(1-\theta'_i)\bm{e}_i,\bm{e})\in  \partial^w\left( \Pexfa \right)$ for all $i=1,\ldots,m$ and 
$(\bm{e},\bm{e}+(\phi'_r-1)\bm{e}_r) \in  \partial^w\left( \Pexfa \right)$ for all $r=1,\ldots,s$.
\end{proof}

\noindent
\textbf{Lemma~\ref{Le3.2}}
\begin{proof}\label{proofLe3.2}
Suppose that  $G^w(\bm{x},\bm{y})=   g(\bm{e},\bm{e}+(\phi'_r-1)\bm{e}_r)$. 
There exists an active constraint $y_r+(\phi'_r-1)y_r=0$ or 
$-\left( \bm{v}^l\right)^{_\top}\bm{x}+\left(\bm{u}^l\right)^{_\top}(\bm{y}+(\phi'_r-1){y_r}\bm{e}_r)-\psi^l=0$ for some $l\in \{1,\ldots,K\}$.  It follows from $y_r>0$ and $\phi'_r\geq 1$ that 
 $y_r+(\phi'_r-1)y_r>0$ and 
 $-\left( \bm{v}^l\right)^{_\top}\bm{x}+\left(\bm{u}^l\right)^{_\top}(\bm{y}+(\phi'_r-1){y_r}\bm{e}_r)-\psi^l=0$ for some  $l\in \{1,\ldots,K\}$.  
 It follows from $(\bm{v}^l,\bm{u}^l)\in \bbR^{m+s}_{++}$ that  $(\bm{x},\bm{y}+(\phi'_r-1){y_r}\bm{e}_r)\in  \partial^s(P_{\scalebox{0.5}{\it EXFA}})$.
 \par
Suppose that  $G^w(\bm{x},\bm{y})=  g(\bm{e}-(1-\theta'_i)x_i\bm{e}_i,\bm{e})$. There exists an active constraint $x_i-(1-\theta'_i)x_i=0$ or 
$-\left(\bm{v}^l\right)^{_\top}(\bm{x}-(1-\theta'_i)x_i\bm{e}_i)+\left(\bm{u}^l\right)^{_\top}\bm{y}-\psi^l=0$ for some $l\in \{1,\ldots,K\}$.  It follows from $x_i>0$ and $\theta'_i\geq 0$ that 
 $\theta'_i=0$ or  $-\left(\bm{v}^l\right)^{_\top}(\bm{x}-(1-\theta'_i)x_i\bm{e}_i)+\left(\bm{u}^l\right)^{_\top}\bm{y}-\psi^l=0$ for some $l\in \{1,\ldots,K\}$. 
 If $\theta'_i=0$, then   it follows from $\phi'_1\geq 1$ that 
 \begin{align*}
 G^w(\bm{x},\bm{y})&= g(\bm{e}-(1-\theta'_i)x_i\bm{e}_i,\bm{e})=\frac{1}{m+s}(m+s-1)=1-\frac{1}{m+s}\\
                         &< \frac{1}{m+s}\left(m+s-1+\frac{1}{\phi'_1}\right) =  g(\bm{e},\bm{e}+(\phi'_1-1)\bm{e}_1)\leq  G^w(\bm{x},\bm{y}), 
 \end{align*}
 which is a contradiction. Therefore, $\theta'_i>0$ and there exists an index $l\in \{1,\ldots,K\}$  such that $-\left(\bm{v}^l\right)^{_\top}(\bm{x}-(1-\theta'_i)x_i\bm{e}_i)+\left(\bm{u}^l\right)^{_\top}\bm{y}-\psi^l=0$.
 It follows from $(\bm{v}^l,\bm{u}^l)\in \bbR^{m+s}_{++}$ that $(\bm{x}-(1-\theta'_i)x_i\bm{e}_i,\bm{y}) \in  \partial^s(P_{\scalebox{0.5}{\it EXFA}})$. 
 \par
Since $\theta'_i=0$ means  $G^w(\bm{x},\bm{y})> g(\bm{e}-(1-\theta'_i)x_i\bm{e}_i,\bm{e})$, 
we have \eqref{eqGw1}  for  any  $(\bm{x},\bm{y})\in P\cap \bbR^{m+s}_{++}$.  
 \end{proof}  

\noindent
\textbf{Theorem~\ref{Th3.3}}
 \begin{proof}\label{proofTh3.3}
It follows from Lemma~\ref{Le3.2} that an optimal solution of the DEA model in~\eqref{maxRMw} is that of the max RM model in~\eqref{p1}--\eqref{p3}, and hence,  
$G^w(\bm{x},\bm{y})= G(\bm{x},\bm{y})$ for any $(\bm{x},\bm{y})\in P\cap \bbR^{m+s}_{++}$. 
\par
Suppose that $(\bm{x}^1,\bm{y}^1), (\bm{x}^2,\bm{y}^2)\in P\cap \bbR^{m+s}_{++}$ satisfying $(\bm{x}^1,-\bm{y}^1)\leq  (\bm{x}^2,-\bm{y}^2)$ and  $(\bm{x}^1,-\bm{y}^1)\not=  (\bm{x}^2,-\bm{y}^2)$. Then, 
$G^w(\bm{x}^2,\bm{y}^2) $ holds for either of the following two cases:  
\begin{align*}
G^w(\bm{x}^2,\bm{y}^2) &= g(\bm{e}-(1-\theta^{2}_{i_2})\bm{e}_{i_2},\bm{e}) & \mbox{ for some } i_2 \in \{ 1,\ldots,m\} \\
G^w(\bm{x}^2,\bm{y}^2) &= g(\bm{e},\bm{e}+(1+\phi^{2}_{r_2})\bm{e}_{r_2}) & \mbox{ for some } r_2 \in \{ 1,\ldots,s\},  
\end{align*}
where
\begin{align*}
\theta^{2}_{i_2} = \min\left\{\, \theta \,\left| \, 
\begin{array}{l}
-\left(\bm{v}^l\right)^{_\top}( \bm{x}^2-(1-\theta)x^2_{i_2}\bm{e}_{i_2})  + \left( \bm{u}^l \right)^{_\top} \bm{y}^2 -\psi^l \leq 0, \, l=1,\ldots,K \\
(1-\theta) x^2_{i_2} \geq 0, \, \theta \geq 0
\end{array}
\,\right\},\right. \\
\intertext{ and }
\phi^{2}_{r_2} = \max\left\{\, \phi \, \left| \, 
\begin{array}{l}
-\left( \bm{v}^l \right)^{_\top} \bm{x}^2+ \left( \bm{u}^l \right)^{_\top}(\bm{y}^2+(\phi-1)y^2_{r_2}\bm{e}_{r_2}) -\psi^{l} \leq 0, \, l=1,\ldots,K \\
(\phi-1)y^2_{r_2}  \geq 0
\end{array}
\,\right\}.\right. 
\end{align*}
It follows from $\theta^{2}_{i_2}>0$ that 
\begin{align*}
(\bm{x}^1-(1-\theta^{2}_{i_2})x^1_{i_2}\bm{e}_{i_2},-\bm{y}^1)\leq  (\bm{x}^2-(1-\theta^{2}_{i_2})x^2_{i_2}\bm{e}_{i_2},-\bm{y}^2)
\intertext{  and } 
(\bm{0},\bm{0})<(\bm{x}^1-(1-\theta^{2}_{i_2})x^1_{i_2}\bm{e}_{i_2},-\bm{y}^1)\not=  (\bm{x}^2-(1-\theta^{2}_{i_2})x^2_{i_2}\bm{e}_{i_2},-\bm{y}^2)
\end{align*}
This means that there exists an index $ l\in \{1,\ldots,K\}$ such that  
\begin{equation}
\begin{array}{r}
0=-\left(\bm{v}^{l}\right)^{_\top}( \bm{x}^2-(1-\theta^{2}_{i_2})x^2_{i_2}\bm{e}_{i_2})+ \left( \bm{u}^{l}\right)^{\top} \bm{y}^2 -\psi^{l} \\
 < -\left(\bm{v}^{l}\right)^{_\top}( \bm{x}^1-(1-\theta^{2}_{i_2})x^1_{i_2}\bm{e}_{i_2})+ \left( \bm{u}^{l}\right)^{\top} \bm{y}^1 -\psi^{l} . 
\end{array} 
\label{eqDelta}
\end{equation}
Let 
\begin{align*}
\theta^1_{i_2} = \min\left\{\, \theta \,\left| \, 
\begin{array}{l}
-\left(\bm{v}^l\right)^{_\top}( \bm{x}^1-(1-\theta)x^1_{i_2}\bm{e}_{i_2})  +\left( \bm{u}^l\right)^{_\top} \bm{y}^1 -\psi^l \leq 0, \, l=1,\ldots,K \\
(1-\theta)x^1_{i_2} \geq 0, \, \theta \geq 0
\end{array}
\,\right\}.\right. 
\end{align*}
Then, it follows from \eqref{eqDelta} that $\theta^{2}_{i_2}<\theta^{1}_{i_2}$. Since $1-\theta^{2}_{i_2}>1-\theta^{1}_{i_2}$, we have 
$g'((1-\theta^{2}_{i_2})\bm{e}_{i_2},\bm{e})<g'((1-\theta^{1}_{i_2})\bm{e}_{i_2},\bm{e})
$, and hence,
\begin{align*}
G(\bm{x}^2,\bm{y}^2) =G^w(\bm{x}^2,\bm{y}^2) &=  g(\bm{e}-(1-\theta^{2}_{i_2})\bm{e}_{i_2},\bm{e})\\ 
                              &<  g(\bm{e}-(1-\theta^{1}_{i_2})\bm{e}_{i_2},\bm{e})\leq G^w(\bm{x}^1,\bm{y}^1)=G(\bm{x}^1,\bm{y}^1). 
\end{align*}
\par
It follows from $\phi^{2}_{r_2}>0$ that 
\begin{align*}
(\bm{x}^1,-\bm{y}^1-(\phi^{2}_{r_2}-1)y^1_{r_2}\bm{e}_{r_2}) \leq  (\bm{x}^2,-\bm{y}^2-(\phi^{2}_{r_2}-1)y^2_{r_2}\bm{e}_{r_2})
\intertext{ and }
(\bm{0},\bm{0})<(\bm{x}^1, \bm{y}^1+(\phi^{2}_{r_2}-1)y^1_{r_2}\bm{e}_{r_2}
)\not=  (\bm{x}^2,\bm{y}^2+(\phi^{2}_{r_2}-1)y^2_{r_2}\bm{e}_{r_2}).
\end{align*}
This means that there exists an index   ${l}\in \{1,\ldots,K\}$ such that   
\begin{equation}
\begin{array}{r}
0=-\left( \bm{v}^l \right)^{_\top} \bm{x}^2+\left(\bm{u}^l \right)^{_\top}(\bm{y}^2+(\phi^{2}_{r_2}-1)y^2_{r_2}\bm{e}_{r_2})-\psi^l \, \\
<-\left( \bm{v}^l \right)^{_\top} \bm{x}^1+\left(\bm{u}^l \right)^{_\top} (\bm{y}^1+(\phi^2_{r_2}-1)y^1_{r_2}\bm{e}_{r_2}) -\psi^{l}.
\end{array} \label{eqDelta2}
\end{equation}
Let 
\begin{align*}
\phi^{1}_{r_2} = \max\left\{\, \phi \,\left| \, 
\begin{array}{l}
-\left(\bm{v}^l\right)^{_\top} \bm{x}^1+ \left( \bm{u}^l\right)^{_\top} ( \bm{y}^1+(\phi-1) y^1_{r_2}\bm{e}_{r_2})  -\psi^{l} \leq 0, \, l=1,\ldots,K \\
(\phi-1)y^1_{r_2} \geq 0, 
\end{array}
\,\right\}.\right. 
\end{align*}
Then, it follows from \eqref{eqDelta2} that $\phi^{2}_{r_2}>\phi^{1}_{r_2}$, and hence,
\begin{align*}
G(\bm{x}^2,\bm{y}^2) =G^w(\bm{x}^2,\bm{y}^2)& = g(\bm{e},\bm{e}+(\phi^{2}_{r_2}-1)\bm{e}_{r_2})\\
&<  g(\bm{e},\bm{e}+(\phi^{1}_{r_2}-1)\bm{e}_{r_2})\leq G^w(\bm{x}^1,\bm{y}^1) =G(\bm{x}^1,\bm{y}^1). 
\end{align*}
\end{proof}

\noindent
\textbf{Theorem~\ref{Th3.4}}
 \begin{proof}
 It follows from $\Pexfa \subseteq \bbR^{m+s}_+$ and $(\bm{x},\bm{y})\in P \cap \bbR^{m+s}_{++}$ that 
 \begin{align*}
 0 &\leq \min\left\{\, \theta \, \left|\, (\bm{x}-(1-\theta)x_i\bm e_i,\bm{y}) \in P_{\scalebox{0.5}{\it EXFA}} \right\} \right. \leq 1 & \mbox{ for all } i=1,\ldots,m \\
 1 &\leq \max\left\{\, \phi\, \left|\, (\bm{x}, \bm{y}+(\phi-1)y_r\bm{e}_r) \in P_{\scalebox{0.5}{\it EXFA}} \right\}. \right. & \mbox{ for all } r=1,\ldots,s
 \end{align*}
 This means that, from  \eqref{thetai} and \eqref{phir}, $\theta'_i=\min\left\{\, \theta \, \left|\, (\bm{x}-(1-\theta)x_i\bm e_i,\bm{y}) \in P_{\scalebox{0.5}{\it EXFA}} \right\} \right.$ and
 $\phi'_r=\max\left\{\, \phi\, \left|\, (\bm{x}, \bm{y}+(\phi-1)y_r\bm{e}_r) \in P_{\scalebox{0.5}{\it EXFA}} \right\} \right.$.
Therefore, it follows from \eqref{D-} and \eqref{D+} that 
 \begin{align}
\max_{i=1,\ldots,m} \theta'_i = 1- D^-(\bm{x},\bm{y})\mbox{ and } \min_{r=1,\ldots,s} \phi'_r = 1+ D^+(\bm{x},\bm{y}). \label{eqB5}
 \end{align} 
 Theorem~\ref{Th3.3}  and~\eqref{eqGw} of  Lemma~\ref{Le3.0} imply from~\eqref{eqB5} that 
 \begin{align*}
G(\bm{x},\bm{y}) &= G^w(\bm{x},\bm{y}) =  
  \max\left\{\,
\begin{array}{l}
\max\left\{\left. g(\bm{e}-{(1-\theta'_i)\bm{e}_i},\bm{e}) \right| \, i=1,\ldots.m \right\} \\ 
\max\left\{\left. g(\bm{e},{\bm{e}+(\phi'_r-1)\bm{e}_r}) \right| \, r=1,\ldots.s \right\}
\end{array} \right\} \\ 
&=  \frac{1}{m+s} \times \max\left\{\,
\begin{array}{l}
\max\left\{\left. m+s-1+ \theta'_i  \right| \, i=1,\ldots.m \right\} \\ 
\max\left\{\left. m+s-1+\frac{1}{\phi'_r}  \right| \, r=1,\ldots.s \right\}
\end{array} \right\} \\
&=  \frac{1}{m+s} \times \max\left\{\,
\begin{array}{l}
m+s-1 + \max\left\{\left. \theta'_i  \right| \, i=1,\ldots.m \right\} \\ 
m+s-1 + \max\left\{\left. \frac{1}{\phi'_r}  \right| \, r=1,\ldots.s \right\}
\end{array} \right\} \\
&=  \frac{1}{m+s} \times \max\left\{\,
\begin{array}{l}
m+s-1 + 1- D^-(\bm{x},\bm{y}) \\ 
m+s-1 + \frac{1}{1+D^+(\bm{x},\bm{y})}
\end{array} \right\} \\
&=\frac{1}{m+s} \times\max\left\{ 
m+s-D^-(\bm{x},\bm{y}),\ m+s-\frac{D^+(\bm{x},\bm{y})}{1+D^+(\bm{x},\bm{y})}
\right\}.
 \end{align*}   
\end{proof}

\noindent
\textbf{Corollary~\ref{Co3.5}}
\begin{proof}
It follows from the Theorem~\ref{Th3.4} that
$\left(\sum_{i=1}^m \theta^{*}_ix_i\bm{e}_i ,\sum_{r=1}^s \phi_r^{*}y_r \bm{e}_r\right)$ is either the input-oriented closest target or the output-oriented one. 
Theorem~\ref{Th3.4} implies  $\bm{\theta}^*\in \bbR^m_{++}$. This means, from $\bm{x}\in \bbR^m_{++}$ that  $\sum_{i=1}^m \theta^{*}_ix_i\bm{e}_i \in \bbR^m_{++}$.
Therefore, it follows from definition~\ref{def1} that  the projection point $\left(\sum_{i=1}^m \theta^{*}_ix_i\bm{e}_i ,\sum_{r=1}^s \phi_r^{*}y_r \bm{e}_r\right)$ is not a free-lunch vector.
\end{proof}

\noindent
\textbf{Theorem~\ref{Th4.2}}
\begin{proof}
For   $(\bm{x},\bm{y})\in P\cap \bbR^{m+s}_{++}$, an optimal solution $(\bm{\theta}^*,\bm{\phi}^*)$ of the max SBM DEA model in~\eqref{maxSBM} with $\theta^*_i=0$ 
satisfies 
\begin{align*}
\left(\sum_{p\not=i} \theta^*_p x_p\bm{e}_p + \theta_i^*(x_i+\tau)\bm{e}_i,\sum_{r=1}^s \phi^*_ry_r\bm{e}_r\right)=
\left(\sum_{p\not=i} \theta^*_p x_p\bm{e}_p, \sum_{r=1}^s \phi^*_ry_r\bm{e}_r\right)\\
=\left(\sum_{p=1}^m \theta^*_px_p\bm{e}_p,\sum_{r=1}^s \phi^*_ry_r\bm{e}_r\right)\in 
\partial^s( P_{\scalebox{0.5}{\it EXFA}})\ \mbox{ for any } \tau > 0.
\end{align*}
This means that $(\bm{\theta}^*,\bm{\phi}^*)$ is a feasible solution of 
\begin{align*}
    \max&  \left\{\, \rho(\bm{\theta},\bm{\phi})   \left| \left(\sum_{p\not=i} \theta_px_p\bm{e}_p+\theta_i(x_i+\tau)\bm{e}_i,  \sum_{r=1}^s \phi_r y_r \bm{e}_r\right) \in \partial^s(P_{\scalebox{0.5}{\it EXFA}}), \eqref{Q3} \right\}\right.,  \label{maxSBMi}  
\end{align*}
for any $\tau>0$, and hence, 
\[
\Gamma(\bm{x},\bm{y})=\rho(\bm{\theta}^*,\bm{\phi}^*)\leq \Gamma(\bm{x}+\tau\bm{e}_i,\bm{y})\ \mbox{ for any } \tau > 0.
\]
\end{proof}

\bibliographystyle{apacite}
\bibliography{Reference}
\end{document}